\documentclass[11pt]{article}

\usepackage{amscd, amsmath, amssymb, fancyhdr, graphicx, epsfig, url, mathrsfs}
\usepackage{graphicx}
\usepackage{epigraph}
\usepackage[backref=page]{hyperref}
\renewcommand*{\backrefalt}[4]{%
	\ifcase #1 (Not cited.)%
	\or        (Cited on page~#2.)%
	\else      (Cited on pages~#2.)%
	\fi}

\hypersetup{
	colorlinks   = true,
	citecolor    = magenta,
	linkcolor    = blue,
	urlcolor     = magenta	
}

\numberwithin{equation}{section}

\newcommand{\version}{version 2.1,\ \ Oct 19, 2026}

\makeatletter
\def\x@arrow{\DOTSB\Relbar}
\def\xlongrightarrowfill@{\arrowfill@\relbar\relbar\longrightarrow}
\newcommand{\xlongrightarrow}[2][]{%
        \ext@arrow 0099\xlongrightarrowfill@{#1}{#2}}
\makeatother

\def\eqref#1{(\ref{#1})}
\newcommand{\goth}{\mathfrak}

\newcommand{\arrow}{{\:\longrightarrow\:}}
\newcommand{\Z}{{\Bbb Z}}
\def\C{{\Bbb C}}
\def\P{{\Bbb P}}

\newcommand{\R}{{\Bbb R}}
\newcommand{\Q}{{\Bbb Q}}
\renewcommand{\H}{{\Bbb H}}
\newcommand{\6}{\partial}
\def\1{\sqrt{-1}\:}

\newcommand{\cntrct}                
{\hspace{2pt}\raisebox{1pt}{\text{$\lrcorner$}}\hspace{2pt}}

\renewcommand{\tilde}{\widetilde}
\renewcommand{\bar}{\overline}
\renewcommand{\phi}{\varphi}
\renewcommand{\epsilon}{\varepsilon}
\renewcommand{\geq}{\geqslant}
\renewcommand{\leq}{\leqslant}

\newcommand{\diam}{{\operatorname{\sf diam}}}
\newcommand{\Teich}{\operatorname{\sf Teich}}

\newcommand{\Comp}{\operatorname{\sf Comp}}

\newcommand{\NS}{\operatorname{\rm NS}}

\newcommand{\Iso}{\operatorname{Iso}}

\newcommand{\Id}{\operatorname{Id}}

\newcommand{\Pic}{\operatorname{Pic}}
\newcommand{\Pos}{\operatorname{Pos}}

\newcommand{\Aut}{\operatorname{Aut}}
\newcommand{\Abs}{\operatorname{Abs}}
\newcommand{\Vis}{\operatorname{Vis}}

\newcommand{\Mon}{\operatorname{\sf Mon}}

\newcommand{\Diff}{\operatorname{\sf Diff}}

\newcommand{\rk}{\operatorname{rk}}

\renewcommand{\Im}{\operatorname{Im}}

\newcommand{\bbZ}{\mathbb{Z}}
\newcommand{\bbQ}{\mathbb{Q}}
\newcommand{\bbR}{\mathbb{R}}
\newcommand{\bbC}{\mathbb{C}}
\newcommand{\bbH}{\mathbb{H}}
\newcommand{\bbP}{\mathbb{P}}
\renewcommand{\AA}{\mathcal{A}}
\newcommand{\CC}{\mathcal{C}}

\newcommand{\OO}{\mathcal{O}}
\newcommand{\PP}{\mathcal{P}}
\newcommand{\MM}{\mathcal{M}}
\newcommand{\KK}{\mathcal{K}}
\newcommand{\MBM}{\mathrm{MBM}}
\newcommand{\cS}{\mathcal{S}}
\newcommand{\st}{\enskip |\enskip}
\newcommand{\wdg}{\wedge}
\newcommand{\BirK}{\mathcal{BK}}
\newcommand{\interior}{{in}}
\newcommand{\PE}{\mathrm{PE}}
\newcommand{\bnd}{\mathscr{B}}
\newcommand{\bndint}{\bnd^{\interior}}
\newcommand{\MCG}{\operatorname{\sf MCG}}

\newcommand{\proof}{\noindent{\bf Proof:\ }}
\newcommand{\pstep}{{\bf Proof. Step 1:\ }}

\newcounter{Mycounter}[section]
\newcounter{lemma}[section]
\renewcommand{\thelemma}{{Lemma \thesection.\arabic{lemma}}}
\newcommand{\lemma}{%
    \setcounter{lemma}{\value{Mycounter}}
    \refstepcounter{lemma}
    \stepcounter{Mycounter}
    {\noindent \bf \thelemma:\ }}

\newcounter{claim}[section]
\renewcommand{\theclaim}{{Claim \thesection.\arabic{claim}}}
\newcommand{\claim}{%
    \setcounter{claim}{\value{Mycounter}}
    \refstepcounter{claim}
    \stepcounter{Mycounter}
    {\noindent \bf \theclaim:\ }}

\newcounter{sublemma}[section]
\newcounter{corollary}[section]
\renewcommand{\thecorollary}{{Corollary \thesection.\arabic{corollary}}}
\newcommand{\corollary}{%
    \setcounter{corollary}{\value{Mycounter}}
    \refstepcounter{corollary}
    \stepcounter{Mycounter}
    {\noindent \bf \thecorollary:\ }}

\newcounter{theorem}[section]
\renewcommand{\thetheorem}{{Theorem \thesection.\arabic{theorem}}}
\newcommand{\theorem}{%
    \setcounter{theorem}{\value{Mycounter}}
    \refstepcounter{theorem}
    \stepcounter{Mycounter}
    {\noindent \bf \thetheorem:\ }}

\newcounter{conjecture}[section]
\newcounter{proposition}[section]
\renewcommand{\theproposition}
      {{Proposition \thesection.\arabic{proposition}}}
\newcommand{\proposition}{%
    \setcounter{proposition}{\value{Mycounter}}
    \refstepcounter{proposition}
    \stepcounter{Mycounter}
    {\noindent \bf \theproposition:\ }}

\newcounter{definition}[section]
\renewcommand{\thedefinition}
      {{Definition~\thesection.\arabic{definition}}}
\newcommand{\definition}{%
    \setcounter{definition}{\value{Mycounter}}
    \refstepcounter{definition}
    \stepcounter{Mycounter}
    {\noindent \bf \thedefinition:\ }}

\newcounter{example}[section]
\renewcommand{\theexample}{{Example \thesection.\arabic{example}}}
\newcommand{\example}{%
    \setcounter{example}{\value{Mycounter}}
    \refstepcounter{example}
    \stepcounter{Mycounter}
    {\noindent \bf \theexample:\ }}

\newcounter{remark}[section]
\renewcommand{\theremark}{{Remark \thesection.\arabic{remark}}}
\newcommand{\remark}{%
    \setcounter{remark}{\value{Mycounter}}
    \refstepcounter{remark}
    \stepcounter{Mycounter}
    {\noindent \bf \theremark:\ }}

\newcounter{problem}[section]
\newcounter{question}[section]
\renewcommand{\thequestion}{{Question \thesection.\arabic{question}}}
\newcommand{\question}{%
    \setcounter{question}{\value{Mycounter}}
    \refstepcounter{question}
    \stepcounter{Mycounter}
    {\noindent \bf \thequestion:\ }}

\makeatletter

\@addtoreset{equation}{section} \@addtoreset{footnote}{section}
\makeatother

\def\blacksquare{\hbox{\vrule width 5pt height 5pt depth 0pt}}
\def\endproof{\blacksquare}

\begin{document}

\begin{center}
{\LARGE\bf
Apollonian carpets and the boundary of \\[2mm]
the K\"ahler cone of a hyperk\"ahler manifold
}

\bigskip

Ekaterina Amerik,\footnote{Supported by 
HSE University basic research program; also partially supported by ANR (France) project FANOHK} Andrey Soldatenkov\footnote{Partially supported by 
FAPESP grant 2024/23819-0 
 }, 
Misha Verbitsky

\end{center}

{\small \hspace{0.05\linewidth}
\begin{minipage}{0.8\linewidth}
\begin{center}
	{\bf Abstract}
\end{center}
The boundary of the ample cone of a hyperk\"ahler
manifold $M$ may contain points which are isotropic 
with respect to the Beauville--Bogomolov--Fujiki form; we call the projectivized set
of such points the isotropic boundary. 
 We are interested in the 
relation between the shape of the isotropic
boundary and the dynamics of the holomorphic
automorphism group $G$ of $M$. The projectivization of the ample cone is 
realized as an open, locally polyhedral subset
in a hyperbolic space ${\mathbb H}$. Up to finite index and quotients by finite subgroups, $G$ is its group of isometries.
The isotropic boundary $S$ is realized as a subset of the hyperbolic
boundary (the absolute) $\operatorname{Abs}$ of ${\mathbb H}$, which is naturally identified
with a Euclidean sphere. It is clear that 
the isotropic boundary $S$ contains the
limit set of $G$ acting on the ample cone.
We prove that, conversely, all irrational
points on $S$ belong to the limit set. 
Using a result of N. Shah about limiting 
distributions of curves under the geodesic 
flow on hyperbolic manifolds, we prove that every real analytic curve in $S$ is contained in a geodesic sphere in $S$,    
and in the presence of such curves the limit set is the closure of the union of these geodesic spheres. We study
the geometry of such configurations of spheres, called Apollonian
carpets, and establish a link
between the geometry of the Apollonian carpet and the
structure of the automorphism group.
\end{minipage}
}



\tableofcontents

\section{Introduction}

{\setlength\epigraphwidth{0.9\linewidth}
\epigraph{
\it \qquad
The abyss and its veil thus came to stand for elements of inner
and outer experience, respectively, both of which were part of the
necessary synthesis. In Belyj's later writings the covering over the
abyss does not merely represent the deception of logical form, but rather
{\bf ``the golden Apollonian carpet''}, which is essential to the creative art
as is the spirit of Dionysus which dwells in the abyss.
}
{\sc\scriptsize Laura Goering, ``The Abyss of Language and
  the Language of the Abyss: Bely's symbolist essays'',
  The Andrej Belyj Society Newsletter, 8, 1989}}


\subsection{The boundary of the ample cone of a hyperk\"ahler
  manifold}


Throughout this paper, we 
treat a hyperk\"ahler manifold
as a complex manifold of K\"ahler type, admitting
a holomorphic symplectic form. Let $M$ be a compact 
hyperk\"ahler manifold of maximal
holonomy (\ref{_hyper_max_holo_Definition_}; see
Subsection \ref{_Basic_Subsection_}
for more details and references).
Then the second cohomology of $M$ is equipped
with a non-degenerate, rational, bilinear
symmetric form $q$ of signature $(3, b_2-3)$,
called {\bf the Beauville--Bogomolov--Fujiki} form,
or the {\bf BBF form} (Subsection \ref{_Basic_Subsection_}).
Complex automorphisms of $M$ act on $H^2(M)$
preserving $q$ and the Hodge decomposition.
Denote by $\KK_M$ the K\"ahler cone, i.e. the set of all K\"ahler classes in $H^{1,1}(M)$, 
and by $\AA_M$ the set of all ample classes in $\NS_\R(M)= (H^{1,1}(M) \cap H^2(M, \Q))\otimes_\Q \R$.
 This set is called {\bf the ample cone of
  $M$}; it is obtained as the intersection of $\NS_\R$
and the K\"ahler cone of $M$.
 
This paper grew out of an attempt to understand
the geometry of the boundary  of $\KK_M$ and $\AA_M$.
These cones were described as follows in 
\cite{_AV:MBM_}. Let $\CC_M \subset  H^{1,1}(M, \R)$
be the connected component of the set of all vectors
with positive square that contains the K\"ahler cone,
and consider the set of so-called ``MBM
classes'' in $H^2(M, \Z)$ (Subsection \ref{sec_cones}). Denote by $R_I$
the set of all effective MBM classes of type (1,1)
(that is, MBM classes represented by a linear
combination of complex curves with positive coefficients, under identification of $H^2(M, \Z)$ with a sublattice in $H_2(M, \Z)$ given by the BBF form). 
Then 
\begin{equation}\label{_Kah_cone_MBM_intro_Equation_}
 \KK_M=\{ \eta\in \CC_M \ \ |\ \ \ \forall u\in R_I, \ \ q(\eta, u) >0
\}
\end{equation}
In other words, the K\"ahler cone is a subset of $\CC_M$
given by a system of linear inequalities. Its boundary
inside $\CC_M$ is locally polyhedral 
(\cite{_AV:Kaw_Mor_}). In this paper, we are
interested in the {\bf isotropic boundary} of
$\KK_M$ and $\AA_M$, that is, the part of the boundary
contained in the set $\{\eta \in H^{1,1}(M, \R)\ \ |\ \ q(\eta,\eta)=0\}.$
We denote the isotropic boundary by $\6\KK_M$ and $\6\AA_M$.

When the Picard number of $M$, denoted by $\rho(M)$, is not maximal,
that is, $\rk( \NS) < b_2-2$, the isotropic boundary $\6\KK_M$
contains an open subset in the cone of isotropic elements $\eta \in H^{1,1}(M)$,
as follows from 
\cite[Proposition 3.4]{_Verbitsky:ergodic_Erratum_}.
However, $\6\AA_M= \6\KK_M\cap \NS_\R$ is a complicated
fractal set, if $R_I$ is non-empty and the Picard number of $M$ is at least 3.
 Further in this paper we relate the
geometry of this fractal with the geometric properties
of the manifold $M$ and its automorphism group.

\subsection{The limit set of the automorphism group}
\label{_limit_intro_Subsection_}

Dynamical properties  
of the group $\Aut(M)$ of holomorphic automorphisms
of a hyperk\"ahler manifold $M$ are closely related to 
the dynamics of the group of isometries of the hyperbolic
space. We always use one of the standard models of the hyperbolic
space, identifying it with the projectivization of
the cone of positive vectors in a quadratic vector space of signature
$(1,n)$. We are mostly interested in the case when
the hyperbolic space is ${\Bbb H}:={\Bbb P}(\CC_M \cap \NS_\R)$, when $M$ is projective of Picard number at least 3.

The action of $\Aut(M)$ on $H^{1,1}(M,\R)$
preserves the BBF form, giving a homomorphism
$\Aut(M)\arrow O(1, b_2-3)$.

Any holomorphic automorphism of $M$
which preserves a K\"ahler class $[\omega]\in \KK_M$, acts on $M$
by isometries of the corresponding Calabi-Yau metric. 
Since the automorphism group of $M$ is discrete,
the group of automorphisms preserving $[\omega]$ is
finite, and $\Aut(M)$ is mapped to $O(1, b_2-3)$
with finite kernel. We are interested in the dynamics of 
its image acting on the corresponding hyperbolic
space ${\Bbb P}\CC_M$, and, even more, on a smaller hyperbolic
space, ${\Bbb P}(\CC_M\cap \NS_\R)$.

Consider a quadratic vector space $(V, q)$ of signature $(1,n)$,
and let ${\Bbb H}:={\Bbb P}(V^+)$ be the corresponding
hyperbolic space.
The {\bf absolute} of ${\Bbb H}$ is the set
$\Abs:= {\Bbb P}(\{\eta \in V\ \ |\ \ q(\eta, \eta)=0\})$,
identified with the Euclidean sphere.
Given a group $\Gamma\subset SO^+(1, n)$ acting on 
${\Bbb  H}$, its {\bf limit set} is the set of all $x\in \Abs$
obtained as limit points of an orbit of $\Gamma$.

It turns out that $\6\AA_M$ is closely related to
the limit set of the automorphism group.
We call a cohomology class
$\eta \in \NS_\R$ {\bf irrational}
if it is not proportional to 
a rational cohomology class.

\hfill

\theorem
Let $M$ be a projective hyperk\"ahler manifold with $\rho(M)\geq 3$,
and $\Aut(M)$ the group of its holomorphic
automorphisms.
Then all irrational points on the 
isotropic boundary $\6 ({\Bbb P}\AA_M)$ 
belong to the limit set of 
$\Aut(M)$ acting on ${\Bbb P}\AA_M$.
Conversely, any point in the limit
set of $\Aut(M)$-action belongs to
$\6 ({\Bbb P}\AA_M)$.

\proof \ref{_Limit_set_Theorem_}. \endproof

\hfill

By \cite[Theorem 1.2.24]{_Cano_Mavarrete_Seade_}, 
all points of the limit set $\Lambda(\Gamma)$
have dense orbits in $\Lambda(\Gamma)$, unless
the group $\Gamma$ is elementary, see \ref{thm_cano}. The 
fractal geometry of $\Lambda(\Gamma)$ is the main subject 
of this paper. 

\subsection{CHOPBs and Apollonian carpets}

The group of automorphisms of a projective hyperk\"ahler
manifold can be described explicitly in terms of 
its MBM classes and Hodge structure, as shown in \cite{_Markman:survey_},
\cite{_AV:Kaw_Mor_} and \cite{_AV:hyperb_}.
Let $\Mon\subset O(H^2(M, \Z), q)$
be the group generated by the monodromy 
of the Gauss-Manin local systems for
all complex deformations of $M$.
This group is called {\bf the monodromy group}
(Section \ref{_monodromy_Section_}).
It has finite index in $O(H^2(M, \Z), q)$
(\cite{_Verbitsky:Torelli_}).
The {\bf Hodge monodromy group}
$\Mon_I$ is its subgroup which
preserves the Hodge decomposition
associated to a complex structure $I$
of hyperk\"ahler type. It is mapped
to $O(H^2(M, \Z)\cap H^{1,1}(M))$ with
finite kernel and finite index for $M$ projective
(\ref{_Mon_I_arithm_Proposition_}).

Let $\Gamma$ be the image of $\Mon_I$ in 
$O(H^2(M, \Z)\cap H^{1,1}(M))$, and $\Gamma_P$ the subgroup of $\Gamma$
which preserves the ample cone. 
As shown by Markman \cite{_Markman:survey_},
the image of the group $\Aut(M)$ of automorphisms of a hyperk\"ahler
manifold in $\Mon_I$ is the subgroup of $\Mon_I$
which preserves the K\"ahler cone.
Therefore the natural map $\Aut(M)\arrow \Gamma_P$ is surjective
with finite kernel, and $\Gamma_P$ is of finite index in the subgroup of isometries of the N\'eron--Severi lattice
$H^2(M, \Z)\cap H^{1,1}(M)$ that preserve the ample cone.  

Let ${\Bbb H}$ be the hyperbolic space obtained
as the projectivization of the positive cone $\Pos$ in $\NS_\R$.
Since $\Gamma$ has finite index in  $O(H^2(M, \Z)\cap H^{1,1}(M))$, 
the quotient ${\cal M}:={\Bbb H}/\Gamma$ is a hyperbolic orbifold
of finite volume. The ample cone can be described as
\[ \AA_M=\{ \eta\in \Pos \ \ |\ \ \ q(\eta, u) >0\ \ \forall u\in R_I\}
\]
as follows from \eqref{_Kah_cone_MBM_intro_Equation_}.
In fact, $\AA_M$ is a connected component of the complement in $\Pos$ to the union of hyperplanes $z^\bot$, where $z$ belongs to the
set $R$ of MBM classes of type $(1,1)$. 
The group $\Gamma$ preserves the set $R$ and
acts on it with finitely many orbits (\cite{_AV:Kaw_Mor_}).
Therefore, the image of the union of the orthogonal complements to $R$
is a finite union of closed hyperbolic orbifolds
${\cal M}_1, \ldots, {\cal M}_k$ immersed
in ${\cal M}={\Bbb H}/\Gamma$ (\ref{_kahler_convex_hy_polyhe_Theorem_}). 
The complement ${\cal M}\backslash\bigcup_i {\cal M}_i$
is a union of finitely many locally convex orbifold 
polyhedral spaces. Each of these orbifolds is obtained
as the quotient of a connected component $U$ of
${\Bbb H}\backslash \bigcup_{z\in R} z^\bot$,
with $U$ being a convex, open subset of ${\Bbb H}$
called {\bf an ample chamber}. Every ample chamber
$P$ is either the ample cone of a birational model of $M$, or
the image of such an ample cone under an element of $\Mon_I$
(\ref{chambers}); in both cases, $P$ is the ample cone of $(M, I_P)$
for an appropriate complex structure $I_P$ on $M$
(\cite{_Markman:survey_}, \cite{_AV:hyperb_}).
The automorphism group of $(M, I_P)$ 
is mapped to the stabilizer of $P$ in $\Mon_I$
with finite kernel and finite index (\ref{_Mon_I_arithm_Proposition_}).

This way, we reduce the study of the automorphism 
group of a hyperk\"ahler manifold to the study of the orbifold fundamental group of
a certain type of a polyhedral space, called CHOPB.

\hfill

\definition
Let ${\cal H}= {\Bbb H}^n/\Gamma$
be a hyperbolic orbifold, and $S_1, \ldots, S_m\subset {\cal H}$
a collection of closed locally geodesic immersed hypersurfaces.
A {\bf convex hyperbolic orbifold with 
 polyhedral boundary} (CHOPB) is 
a connected component of the complement
${\cal H}\backslash \bigcup_i S_i$. 

\hfill

Let ${\cal P}$ be a CHOPB, and $P\subset {\Bbb H}^n$
a connected component of its preimage in ${\Bbb H}^n$.
{\bf The isotropic boundary} of the CHOPB
is the set $\bar P\cap \Abs$, where 
$\Abs$ is the absolute of ${\Bbb H}^n$. It turns out that the
isotropic boundary of a CHOPB has empty interior in $\Abs$,
unless $P={\Bbb H}^n$ (cf. \ref{_Limit_set_Theorem_}).

As mentioned in Subsection
\ref{_limit_intro_Subsection_}, 
the limit set $\Lambda$ of $\Gamma_P= \pi_1({\cal P})$
is in the isotropic boundary of $P$, and,
conversely, all non-cusp points of the
isotropic boundary belong to $\Lambda$
(\ref{_Limit_set_Theorem_}). 

The nature of this set is explored in 
some detail in Section \ref{_Ap_carpet_Section_}.
We define {\bf the Apollonian carpet} of 
the isotropic boundary of a CHOPB as 
the union of all positive-dimensional connected
real analytic subvarieties in its isotropic boundary,
and prove that each component of the 
Apollonian carpet is a geodesic sphere $S$ in $\Abs$
(\ref{_carpet_union_Theorem_}). Moreover, the
convex hull of each $S$ is a hyperbolic subspace 
in the closure of $P$, and its image in the quotient
${\Bbb H}^n/\Gamma$ is a closed
hyperbolic immersed submanifold ${\cal N}$ in  the
closure of ${\cal P}$.
This implies, in particular, that the lattice
$\pi_1({\cal N})$ is a subgroup of $\pi_1({\cal P})$
(\ref{_Apollonian_carpet_explicitly_Theorem_}; here
this result is stated for CHOPBs associated with the
hyperk\"ahler manifolds, but the
proof works generally).

\hfill

\remark\label{_thin_groups_Remark_}
The results of Sections \ref{chamb}, \ref{_limit_points_Section_} and
\ref{_carpet_Section_} do not use hyperk\"ahler geometry.
They apply to any group $\Gamma_P$ obtained as follows:
$\Gamma$ is a subgroup of $SO^+(V)$ commensurable with $SO(V_\Z)$, where
$V_\Z$ is a lattice with a rational quadratic form of signature $(1,n)$,
$W_1,\dots,W_m\subset V$ are rational hyperplanes, and $\Gamma_P$ is the
stabilizer in $\Gamma$ of a connected component $P$ of
${\Bbb H}^n\backslash\bigcup_i\Gamma\cdot{\Bbb P}W_i$; in other words,
$\Gamma_P$ is the orbifold fundamental group of a CHOPB in an arithmetic
hyperbolic orbifold. Unless $P={\Bbb H}^n$, the group $\Gamma_P$ has
infinite index in $\Gamma$, and when its limit set is not contained in a proper
subsphere of $\Abs$ it is Zariski dense in $SO(V)$; that is, $\Gamma_P$ is a thin
group in the sense of Sarnak. This class contains the classical Apollonian group
(Subsection \ref{_example_Apollonian_Subsubsection_}) and is closely related to the
symmetry groups of crystallographic and Kleinian sphere packings
(\cite{_Kontorovich_Nakamura_}, \cite{_Kapovich_Kontorovich_}); however, not
every thin group belongs to this class.
Since the authors are mostly interested in the geometry of hyperk\"ahler manifolds,
we state our results in the setting when $\Gamma_P$ is the group of automorphisms
of such a manifold.

\subsection{Indra's pearls, Apollonian gaskets and
 geometry of K3 surfaces (after A. Baragar et al.)}

The ``Apollonian gasket'' is a relatively modern concept which has its roots in antiquity with
 Apollonius of Perga (c. 240 BC -- c. 190 BC),
see \cite{_Indra:MSW_}. 
It is obtained by taking 
three kissing circles, inscribing a fourth circle
kissing these three, then inscribing 
a circle inside a triangle formed by the
segments of three circles, and so on,
ad infinitum:\\[1mm]

\centerline{\epsfig{file=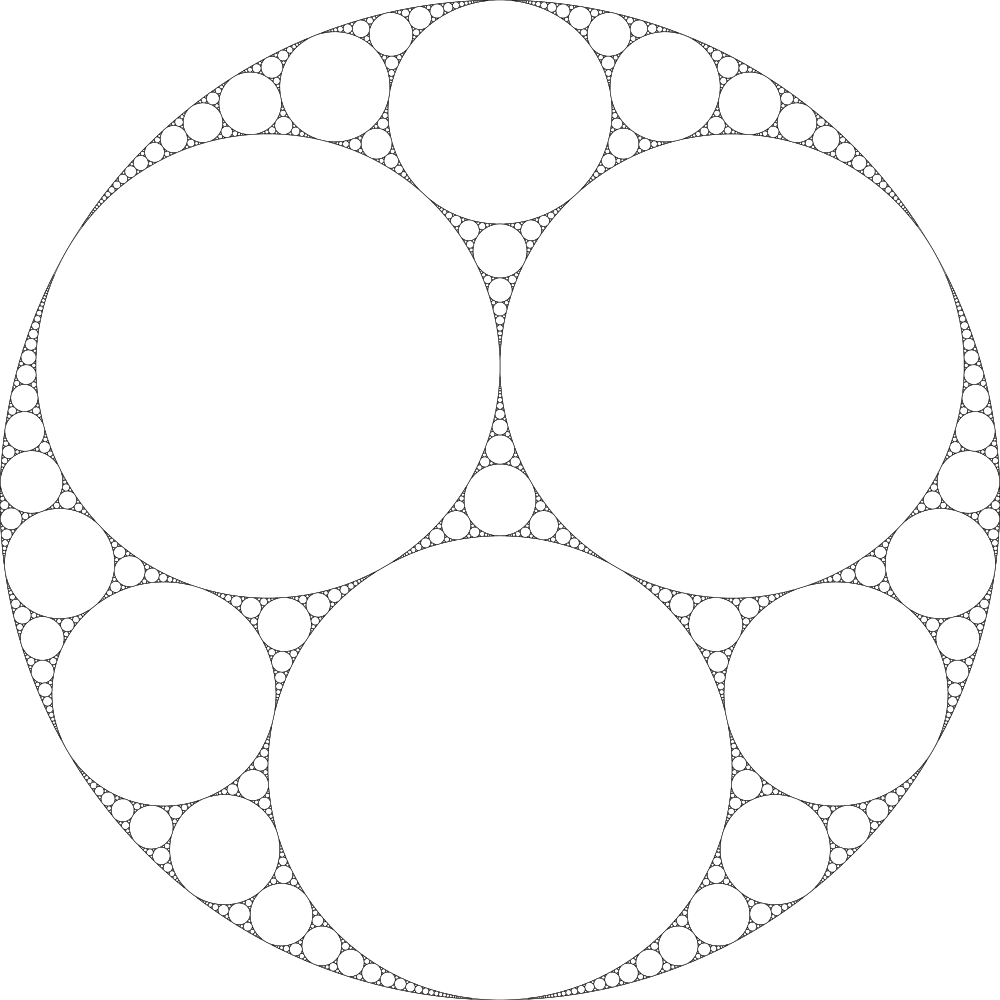,width=0.4\linewidth}}

The resulting configuration of circles is called the Apollonian circle
packing, and its closure is a fractal known as the
Apollonian gasket, see e.g. \cite{_Kontorovich_Nakamura_}, \cite{St}, \cite{VW}.
This construction can (and should) be performed on a 2-sphere.\footnote{A referee points out that especially when related to the geometry of K3 surfaces, some aspects are best seen on the plane,
see e.g. Baragar's works in the references.}
Since the M\"obius group acts transitively on triples of
pairwise kissing circles, the Apollonian gasket on $\C P^1$ 
is unique up to conformal equivalence, see \cite[p. 201]{_Indra:MSW_}.

Interestingly enough, the Apollonian gasket is the limit
set of a discrete subgroup of the M\"obius group $PSL(2, \C)$
acting on $\C P^1$. This is why the Apollonian gasket is one
of the central subjects of ``Indra's Pearls'' \cite{_Indra:MSW_}, by
D. Mumford, C. Series and D. Wright. In \cite{_Indra:MSW_} this was
shown by an explicit argument with matrices, and it also follows
from our \ref{_carpet_union_Theorem_} together with \ref{_Limit_set_Theorem_}.

The Apollonian gasket is ubiquitous in physics
(\cite{_Satija_}) and even in art: the largest
artwork in history, ``Earth Drawings'' by Jim Denevan,
drawn in the sand of Black Rock Desert (Nevada),
depicts the Apollonian circle packing.\\[1mm]

\centerline{\epsfig{file=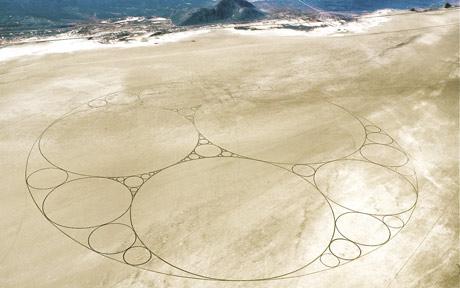,width=0.55\linewidth}}

The Apollonian gasket can be obtained as the isotropic boundary
of a universal covering of a CHOPB.
Consider a circle $S$ in
$\C P^1=\Abs {\Bbb H}^3$. The convex hull of $S$ in ${\Bbb H}^3$
is a hyperbolic 2-plane, denoted by $W$, and $S$ is its isotropic
boundary. The 2-plane $W$ can be obtained as a projectivization
of the positive cone in a subspace $\tilde W= \R^{1,2}$ in $\R^{1,3}$.
Choose an integer lattice in $\R^{1,3}$, and
consider the corresponding lattice group $\Gamma:=SO_\Z(1,3)$.
If $\tilde W$ is rational,
the stabilizer $\Gamma_W$ of $S$ in $\Gamma$ is 
a group which is commensurable with $SO_\Z(1,2)$.
We can choose $\tilde W$ in such a way that
its images under $\Gamma$ can be tangent, but 
never intersect transversally (Subsection
\ref{_example_Apollonian_Subsubsection_}).\footnote{A similar example
was obtained earlier in \cite{_Baragar1_}.}
Then its image in the hyperbolic orbifold ${\cal M}={\Bbb H}^3/\Gamma$
is ${\Bbb H}^2 /\Gamma_W$, a hyperbolic orbifold of finite volume. The
universal cover of the corresponding CHOPB ${\cal P}$ is  a locally finite
hyperbolic polyhedron, which looks as follows:

\centerline{\epsfig{file=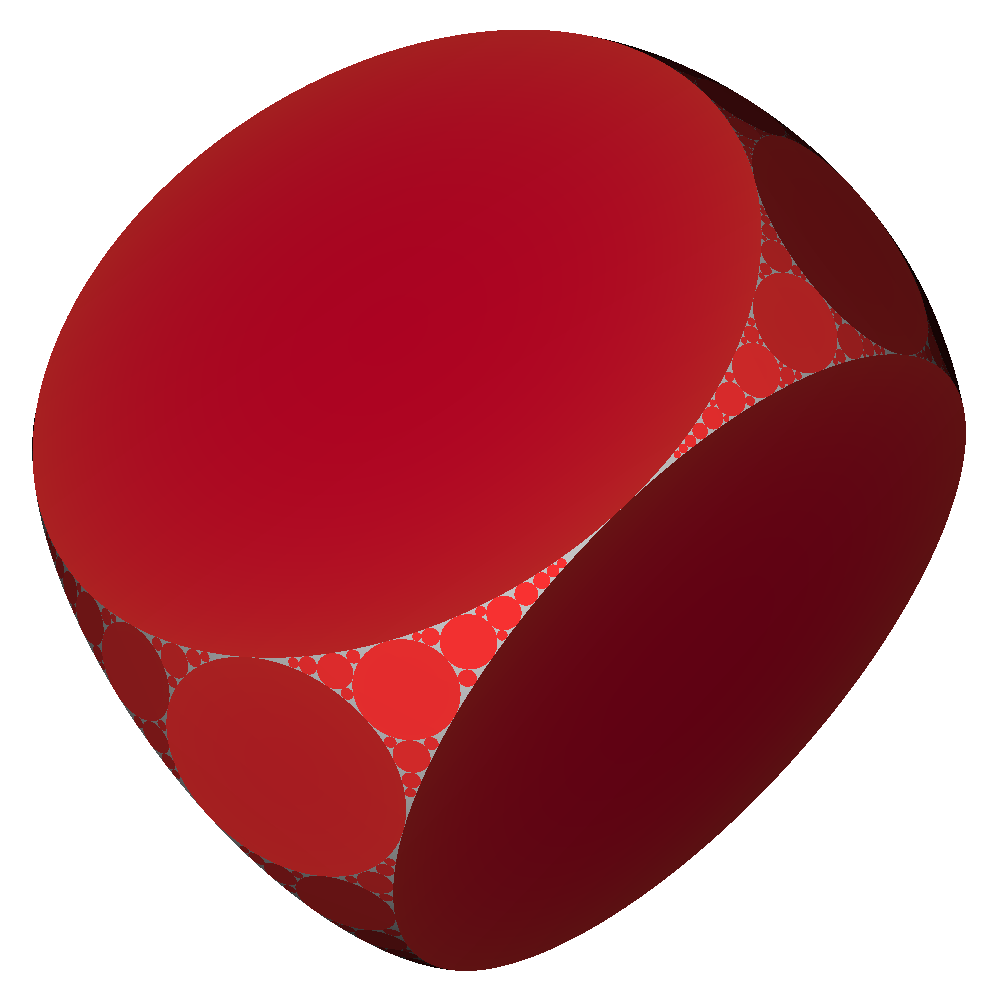,width=0.55\linewidth}}

\noindent
(for a clearer picture, here
we used the Klein model of a hyperbolic space;
 in the Klein model, totally geodesic hypersurfaces
are represented by plane discs).

From this picture, it is visually clear that
the isotropic boundary of ${\cal P}$ is
the Apollonian gasket. A formal proof of this
observation can be obtained by applying
\ref{_Apollonian_carpet_explicitly_Theorem_},
which describes the Apollonian carpet
of a CHOPB in an arithmetic hyperbolic orbifold
(the statement is made for CHOPBs of hyperk\"ahler
origin, but the proof works for this general
situation).

This construction was realized in an ample
cone of a K3 surface by A. Baragar in \cite{_Baragar:Apollonian_}.
In a series of papers (\cite{_Baragar1_,_Baragar:higher_,_Baragar:Apollonian_,_Baragar:7-8_,_Baragar:Enriquez_}), Baragar constructed many examples
illustrating the fractal nature of the isotropic
boundary of the ample cone of a K3 surface.

For K3 surfaces, an MBM class is an integer
class of square $-2$, usually abbreviated as ``$(-2)$-class''.
The BBF form $q$ on a K3 surface is the intersection form.
The corresponding hyperbolic orbifold is
${\Bbb P}(\Pos \cap \NS_\R)/\Gamma$, where $\Gamma$ is a subgroup
of finite index in $SO(\NS_\Z, q)$, and the corresponding
hyperbolic polyhedron $P\subset {\Bbb P}(\Pos\cap\NS_\R)$
is one of the connected components of 
${\Bbb P}(\Pos\cap\NS_\R)\backslash \bigcup_i z_i^\bot$,
where $\{z_i\}\subset \NS_\Z$ is the set of all $(-2)$-classes.

In \cite{_Baragar:Apollonian_}, Baragar constructed an even
lattice $\Lambda$ of rank 4, such that its isometry group acts
on the set of $(-2)$-classes transitively. This means that
$O(\Lambda)$ acts on the set of faces of $P$ transitively.
Moreover, for any pair of $(-2)$-classes  $a, b\in \Lambda$, 
the lattice $\langle a, b \rangle$ generated by $a$ and $b$
has signature $(1,1)$ (this means that their orthogonal
complements don't intersect in the absolute) or $(0,1)$ (and this means
that the orthogonal circles are tangent). This gives a
fractal formed by kissing circles, which is conformally
equivalent to the Apollonian gasket by 
\cite[p. 201]{_Indra:MSW_}, as we indicated earlier.

In his earliest work on the subject 
(\cite{_Baragar1_}), Baragar established the
fractal nature of the isotropic boundary of the ample
cone of a K3 surface, and computed the Hausdorff dimension
of this fractal, identified with the limit set of the
automorphism group. This computation is based on 
a result of D. Sullivan, relating the Hausdorff 
dimension of the limit set with the growth rate
of the number of orbit points in a ball of radius $r$, see \cite{Su}.
See also \cite{VW} and references therein for what
is known about the Hausdorff dimension of the classical
Apollonian gasket.

This observation was later used by I. Dolgachev
\cite{_Dolgachev_}, who used the bounds on the Hausdorff dimension
to estimate the growth rate of the number of curves on a K3 surface
in a given orbit of the automorphism group.
Note that the framework built in the present paper can be used
to extend these estimates to an arbitrary hyperk\"ahler 
manifold. Indeed, the Hausdorff dimension of the isotropic boundary of an ample
cone is bounded from below by the dimension of the
components of its Apollonian carpet.

\subsection{The structure of this paper}

The paper is organized as follows. In Section \ref{HK} we
recall some basic definitions and facts about the
hyperk\"ahler manifolds. We also introduce the convex 
cones in $H^{1,1}(M,\R)$ which are central for our 
further investigation. In Section
\ref{chamb} we prove that whenever the ample, movable or
effective cone does not coincide with the
positive cone in the real N\'eron--Severi group, its isotropic
boundary is nowhere dense (the cone has no ``round
part''); on the contrary, when the Picard number is not
maximal, the K\"ahler cone always has a ``round
part'', which is open in the boundary of the positive 
cone. The statement for the ample cone is mentioned in
\cite{Mat}, with a partial proof, and for the effective
cone it is proved in \cite{_Denisi_}. Our starting point was an
earlier version of \cite{_Denisi_} where the proof for the
effective cone was given for the known families of hyperk\"ahler
manifolds only, by adapting  Kov\'acs's argument \cite{Ko1}
on K3 surfaces. It seems that the simple general proof of this result
relying only on hyperbolic geometry was not widely known. We
have therefore decided to give a self-contained and
complete account here. 

In Sections \ref{_monodromy_Section_} and \ref{_limit_points_Section_}, we study the boundary of the ample
cone using hyperbolic geometry. In Section \ref{_carpet_Section_}, we recall a
theorem of N. Shah which permits us to prove that the
``analytic'' part of the absolute of the ample cone is a
union of spheres. In Section \ref{_Ap_carpet_Section_}, we describe these spheres
more explicitly and give a number of examples of their possible
behavior. We also make a connection to the structure of
the automorphism group of the hyperk\"ahler manifold. In Appendix \ref{app_symmetries},
we describe the symmetry groups of some of the examples.



\section{Hyperk\"ahler manifolds}\label{HK}

We briefly recall the basic definitions of hyperk\"ahler geometry.
We refer to \cite{_Beauville_}, \cite{Hu} and \cite{_Markman:survey_} for a more detailed
discussion.

\subsection{The basic definitions}
\label{_Basic_Subsection_}

Let $M$ be a complex manifold with a Riemannian metric $g$. Denote by $I$ the complex structure on $M$,
and recall that $g$ is Hermitian if for any pair of tangent vectors $u$ and $v$ we
have $g(u,v) = g(Iu,Iv)$. In this case one defines $\omega(u,v) = g(Iu,v)$. Then
$\omega$ is a two-form on $M$, and $M$ is called K\"ahler if $d\omega = 0$.
If $M$ is K\"ahler then $\omega$ is a symplectic form on $M$.

\hfill

\definition
A {\bf hyperk\"ahler structure} on a manifold $M$ is a Riemannian metric $g$
and a triple of complex structures $I$, $J$, $K$, satisfying quaternionic relations
$IJ = - JI =K$, such that $g$ is K\"ahler for $I$, $J$ and $K$.
A {\bf hyperk\"ahler manifold} is a manifold with a hyperk\"ahler structure.

\hfill

A hyperk\"ahler manifold carries symplectic forms $\omega_I$, $\omega_J$ and $\omega_K$
associated to $g$ and the complex structures $I$, $J$ and $K$, as above.
Hyperk\"ahler manifolds are holomorphically symplectic. Indeed, $\Omega_I = \omega_J+\1\omega_K$
is a holomorphic symplectic form on $(M,I)$. The holomorphic symplectic
forms on $(M,J)$ and $(M,K)$ are defined analogously.

Recall that for a Riemannian metric $g$ and a point $x\in M$
the subgroup of $GL(T_x M)$ generated by the operators of parallel
transport (along all loops with endpoints at $x$, with respect to the
Levi--Civita connection) is called the holonomy group of $g$.
A hyperk\"ahler manifold can be equivalently defined as a Riemannian manifold
whose holonomy is contained in $Sp(n)$ (the group of all orthogonal isometries
of a quaternionic vector space preserving $I,J,K$).

In this paper we will consider only one special class of hyperk\"ahler
manifolds.

\hfill

\definition\label{_hyper_max_holo_Definition_}
A compact hyperk\"ahler manifold $M$ is {\bf of maximal holo\-nomy}, or {\bf IHS},
if $\pi_1(M)=1$ and the holonomy group of the Levi--Civita connection
equals $Sp(n)$, where $2n = \dim_\bbC (M)$.

\hfill

For the rest of this paper, we abuse the language
slightly, and ``a hyperk\"ahler manifold'' for us always
means a {\bf compact hyperk\"ahler manifold of maximal holonomy.}

On any hyperk\"ahler manifold of maximal holonomy, the space
of holomorphic two-forms is one-dimensional, spanned by the holomorphic
symplectic form described above.

Fix a hyperk\"ahler manifold $M$ and consider its singular cohomology
group $H^2(M,\bbQ)$. One of the main properties of hyperk\"ahler manifolds
of maximal holonomy is the existence of a canonical quadratic form
$q\in \mathrm{Sym}^2H^2(M,\bbQ)^*$, called the {\bf Beauville--Bogomolov--Fujiki} form
or the {\bf BBF form}. This form can be characterized by the following
property (see \cite{So} for a more detailed discussion):

\hfill

\theorem (Fujiki, \cite{_Fujiki:HK_})
For any hyperk\"ahler manifold $M$ of complex dimension $2n$ there exists
a constant $c_M\in\bbQ$ and a quadratic form $q$ on $H^2(M,\bbQ)$ such that for any
$a\in H^2(M,\bbQ)$ we have
$$\int_M a^{2n}=c_M q(a)^n.$$

\hfill

We may assume, after rescaling, that the BBF form $q$ is integral and primitive on $H^2(M,\Z)$.
In the above theorem the integral on the left-hand side denotes the pairing
of the cohomology class $a^{2n}\in H^{4n}(M,\bbQ)$ with the fundamental
class of $M$ determined by the canonical orientation of $M$, see \cite{So} for details.

The form $q$ has signature $(3,b_2(M)-3)$. If $I$ is one of the complex structures that appears
as part of a hyperk\"ahler structure on $M$, then $q$ is negative definite on the primitive forms
of $I$-type $(1,1)$ and positive definite on the subspace of $H^2(M,\bbR)$
spanned by $\omega_I$, $\mathrm{Re}(\Omega_I)$ and $\mathrm{Im}(\Omega_I)$,
where $\omega_I$ is the K\"ahler form and $\Omega_I$ is the holomorphic symplectic form.

\subsection{Cones in the cohomology of a hyperk\"ahler manifold}\label{sec_cones}

Let us fix a hyperk\"ahler manifold $M$ of maximal holonomy with a
hyperk\"ahler structure $(g,I,J,K)$. We will consider the Hodge decomposition
of the cohomology of $M$ with respect to $I$ if not stated otherwise.

The main goal of the present paper is to study the shape of several natural
convex cones in the real vector space $H^{1,1}(M,\bbR)$. The latter
means the intersection of $H^{1,1}(M)$ with $H^2(M,\bbR)$ inside $H^2(M,\bbC)$.
As was explained above, the vector space $H^{1,1}(M,\bbR)$ carries the quadratic form $q$,
the restriction of the BBF form to the $(1,1)$-part of the cohomology.
The form $q$ has signature $(1,b_2(M) - 3)$ on $H^{1,1}(M,\bbR)$.

Recall that the N\'eron--Severi group $\NS(M)$ is the image of the Picard group $\Pic(M)$
in $H^2(M,\bbZ)$ under the map that computes the first Chern class of a line bundle.
By Hodge theory, the group $\NS(M)$ is contained in $H^{1,1}(M,\bbR)$, and we will denote
by $\NS_\bbR$ the $\bbR$-subspace of $H^{1,1}(M,\bbR)$
spanned by $\NS(M)$.
The condition that a hyperk\"ahler manifold $M$ is projective
is equivalent to the condition that $q|_{\NS_\bbR}$ has signature $(1,\rho(M)-1)$,
where $\rho(M)$ is the Picard number of $M$, see \cite{Hu}.

Let us recall the definitions of the cones that we are going to study.

\begin{enumerate}

\item The positive cone $\CC_M\subset H^{1,1}(M,\bbR)$ is the connected component of the cone
$\{v\in H^{1,1}(M, \bbR)\st q(v)>0\}$ that contains the class of the K\"ahler form $\omega_I$.

\item The K\"ahler cone $\KK_M\subset \CC_M$ is the cone consisting of the cohomology classes of all K\"ahler forms on $M$.
The cone $\KK_M$ is open and its closure is called the nef cone.

\item The birational K\"ahler cone $\BirK_M \subset \CC_M$ is the cone consisting of the cohomology
classes of the form $f^*[\omega]$, where $f\colon M\dashrightarrow M'$ is an arbitrary bimeromorphic map
to a hyperk\"ahler manifold $M'$ and $\omega$ is an arbitrary K\"ahler form on $M'$. Note that
		the cone $\BirK_M$ is not convex (but its closure is, as we shall see from its description below).

\item The pseudoeffective cone $\PP_M\subset H^{1,1}(M, \bbR)$ consists of the cohomology classes
of all closed positive $(1,1)$-currents on $M$. It is known that $\PP_M$ is closed and contains $\CC_M$.

\item The ample cone $\AA_M\subset \NS_\bbR(M)$ is spanned by the classes of ample line bundles.
The ample cone is open in $\NS_\bbR(M)$.

\item The movable cone $\MM_M\subset \NS_\bbR(M)$ is spanned by the classes of movable line bundles, i.e.
the line bundles whose complete linear systems have base loci of codimension at least two.


\end{enumerate}


We summarize the relations between the cones introduced above. We have the following inclusions:
$$
\KK_M\subset\BirK_M\subset \CC_M\subset\PP_M.
$$
Moreover, it is known that $\overline{\BirK_M}$ is dual to $\PP_M$ with respect to $q$, see \cite[Corollary 4.6]{HuCone}.

The ample and the movable cones are related to the K\"ahler cone and the birational K\"ahler cone as follows.
For the ample cone we have $\AA_M = \KK_M \cap \NS_\bbR(M)$.
For the movable cone $\MM_M^\interior = (\overline{\BirK_M})^\interior\cap \NS_\bbR(M)$,
where the superscript $\interior$ means the interior of the cone, see \cite[Proposition 5.6 and Lemma 6.22]{_Markman:survey_}.

The K\"ahler cone and the birational K\"ahler cone admit an important description
in terms of wall and chamber decomposition of the positive cone. To give this description
we need to recall the notions of an MBM class and of a prime exceptional divisor.

A cohomology class $z\in H^2(M,\Z), \ q(z)<0$ on a
hyperk\"ahler manifold $M$ with a complex
structure $I$ is called {\bf an MBM class}\footnote{acronym for ''monodromy birationally minimal''} if some
(equivalently, any) deformation $(M,I')$ of $(M,I)$
with $\NS_\R$ generated by $z$ contains a rational curve;
 see \cite{_AV:survey_}.
We will denote by $\MBM(M)\subset H^2(M,\bbZ)$ the set of primitive MBM classes. The set $\MBM(M)$ is a deformation invariant of
a hyperk\"ahler manifold $M$.

\hfill

\example
Recall that a $(-2)$-class on a K3 surface
$M$ is a cohomology class $\eta\in H^2(M, \Z)$ with
self-intersection $-2$. The primitive MBM classes on a K3 surface
are exactly the $(-2)$-classes (\cite{_AV:MBM_}).

\hfill

Recall from \cite{_AV:Kaw_Mor_} that the set of primitive MBM classes on a hyperk\"ahler
manifold $M$ is $q$-bounded: there exists a constant $\alpha>0$ such that
for all $v\in \MBM(M)$ we have $-\alpha \le q(v) < 0$.
Now consider the intersection $\cS = \MBM(M)\cap \NS(M)$. Then
the collection of hyperplanes $v^\perp$, $v\in \cS$ is locally
finite in $\CC_M$ (see e.g. \cite[Proposition 5.14]{SSV}) and
the complement $\CC_M\setminus \cup_{v\in \cS} v^\perp$ is a disjoint
union of open chambers, the K\"ahler cone $\KK_M$ being one of
the chambers, see \cite{_AV:MBM_}.
In other words, the MBM classes are orthogonal to the
faces of the K\"ahler chambers, and all faces of the K\"ahler
cones of deformations and birational models of $M$ are 
realized in this way (\cite{_AV:hyperb_}).

\hfill

Now we recall an analogous description of the closure of the birational
K\"ahler cone $\overline{\BirK_M}$. Recall that a {\bf prime exceptional
divisor} is a prime divisor $D$ on $M$ whose cohomology class is negative:
$q([D])<0$. Prime exceptional divisors are uniruled, see \cite{Bou}. We will denote by $\PE(M)\subset \NS(M)$
the set of cohomology classes of prime exceptional divisors on $M$.

It follows from the uniruledness (see \cite[Lemma 2.11]{_AV:survey_})
that the cohomology class $[D]$ of a prime exceptional divisor $D$
is MBM. We will call an MBM class of type $(1,1)$ {\bf divisorial}
if its orbit under the action of the Hodge monodromy group (see \cite{_Markman:survey_})
contains a class that is proportional to the class of a prime exceptional divisor.
Let us now denote by $\cS$ the set of divisorial primitive MBM classes on a
hyperk\"ahler manifold $M$. The complement $\CC_M\setminus \cup_{v\in \cS} v^\perp$
is again a disjoint union of open chambers and the closure of one
of these chambers is $\overline{\BirK_M}$.

Let us finally recall the following description of the pseudoeffective
cone given in \cite{Bou}:
\begin{equation}\label{eqn_pseudoeffective}
\PP_M = \overline{\CC_M + \sum_{v\in \PE(M)}\bbR_+\cdot v},
\end{equation}
which means that the boundary of the cone $\PP_M$ is locally polyhedral away from $\CC_M$,
with extremal rays spanned by the classes of the prime exceptional
divisors.


\section{Hyperbolic spaces and arithmetic lattices}\label{chamb}


In this section we recall some standard notions of hyperbolic geometry
in the form convenient for applications to hyperk\"ahler manifolds.
As a general reference on hyperbolic geometry one may use, for example,
\cite{_Kapovich:hyperbolic_}, where most of the facts that we recall
are proven in greater generality. In order to streamline the exposition,
we choose not to deduce what we need from those general theorems,
but to give self-contained elementary proofs, claiming no originality
whatsoever. We also immediately apply the statements about the hyperbolic
spaces to hyperk\"ahler manifolds, reformulating them in the terminology
introduced in the previous section.

\subsection{The hyperbolic space}\label{sec_hyperbolic}

Let $V$ be a vector space over $\bbR$ and $q\in \mathrm{Sym}^2V^*$
a quadratic form of signature $(1,d)$.
Let $V^+ = \{x\in V\st q(x)>0\}\subset V$ be the cone of positive vectors 
and $\bbP (V^+)$ its image in the real projective space $\bbP (V)$.
Analogously, set $V^- = \{v\in V\st q(v)<0\}$.
Denote by $Q = \{x\in \bbP (V)\st q(x)=0\} \subset \bbP (V)$ the
real projective quadric defined by the quadratic form $q$.
The quadric $Q$ is the boundary sphere of the hyperbolic space $\bbP(V^+)$.
For a non-zero vector $v\in V$ we will denote by $[v]$ its image in $\bbP(V)$.

The group $O(V,q)$ has 4 connected components and
we denote the connected component of the identity by $SO^+(V,q)$.
It is clear that $\bbP(V^+)\simeq SO^+(V,q)/SO(d)$.
The hyperbolic space $\bbP(V^+)$ carries an $SO^+(V,q)$-inva\-riant
Riemannian metric $g$. It is easy to see that this metric
is unique up to a constant multiplier. One can also check
that $SO^+(V,q)$ is the group of all oriented isometries of
$\bbH^d=(\bbP(V^+), g)$. We recall the following classification
of isometries of a hyperbolic space.

\hfill

\theorem\label{_isometries_classification_Theorem_}
(Classification of isometries; this theorem also serves as a definition.)
Let $\alpha \in SO^+(1,d)$ be an isometry acting on $\bbH^d$.
Then one and only one of the following cases occurs:
\begin{enumerate}
\item $\alpha$ has an eigenvector $x$ with $q(x)>0$ (then $\alpha$ is called {\bf elliptic});
\item $\alpha$ has an eigenvector $x$ with $q(x)=0$ and eigenvalue $\lambda$ satisfying $|\lambda|>1$
(then $\alpha$ is called {\bf hyperbolic}, or {\bf loxodromic});
\item $\alpha$ has a unique eigenvector $x$ with $q(x)=0$ (then $\alpha$ is called {\bf parabolic}).
\end{enumerate}
%
%
%
%

\subsection{Wall and chamber decompositions of the hyperbolic space}

For the applications to hyperk\"ahler geometry, we need to study wall and chamber
decompositions of the hyperbolic space $\bbP(V^+)$. Such decompositions are
determined by a subset $\cS\subset V^-$ of elements orthogonal to the walls.
We will always assume the following:
\begin{equation}\label{assumptions_S}
\cS\neq\emptyset,\quad -\cS=\cS,\quad \cS \mbox{ is discrete in } V \mbox{ and $q$-bounded}.
\end{equation}
The latter condition means that there exists some constant $\alpha>0$ such that $-\alpha \le q(v) < 0$ for any $v\in \cS$.

For $v\in V^-$ denote $H_v = v^\perp\subset V$ and
\begin{equation}\label{eqn_vdiamond}
V^\Diamond = V^+\setminus \cup_{v\in \cS} H_v.
\end{equation}
It is well-known that the collection of hyperplanes $H_v$, $v\in \cS$ is locally finite in $V^+$,
see e.g. \cite[Proposition 5.14]{SSV}. Therefore $V^\Diamond$ is a union of disjoint open
connected components which we will call the {\bf $\cS$-chambers} of $V^+$. 

Given an $\cS$-chamber $\CC$ we will denote by $\bbP(\CC)$ its image in $\bbP(V)$.
We are interested in the structure of the boundary of $\bbP(\CC)$ at infinity,
or, more precisely, in the structure of the set $\bbP(\overline{\CC})\cap Q$.
The properties of this set depend on the geometry of $\cS$, and there are
two substantially different cases to consider: either the set $\cS$ is
``big'', i.e. spans $V$ and is invariant under the action of an arithmetic
lattice in $SO(V,q)$, or $\cS$ is contained in a proper subspace $N\subset V$.
In the applications, the first case arises when $V$ is the space spanned by the
N\'eron--Severi lattice and $\CC$ is the ample or the movable cone, and the second case arises
when $V = H^{1,1}(M,\bbR)$ and $\CC$ is the K\"ahler cone (the subspace
$N$ is then spanned by the N\'eron--Severi lattice). We will treat these
two cases separately below.

\subsubsection{The case when $\cS$ spans $V$}

In the above setting, assume moreover that $V = V_\bbZ\otimes \bbR$ for some
lattice $V_\bbZ$. Denote $V_\bbQ = V_\bbZ\otimes \bbQ$ and assume that $q$
is rational, i.e. $q\in \mathrm{Sym}^2V_\bbQ^*$. Let $\Gamma\subset O(V_\bbZ,q)$ be
a subgroup of finite index and assume that $\cS\subset V_\bbZ$ and 
$\cS$ is $\Gamma$-invariant.
The following proposition is well known.
We sketch an elementary proof for the sake of completeness.

\hfill

\proposition\label{prop_dense}
In the above setting assume that $d\ge 2$.
For any $x\in Q$ the orbit $\Gamma x$ is dense in $Q$ in the analytic topology.

\hfill

\begin{proof}
Let $U\subset Q$ be a non-empty open subset. We have to prove that $\Gamma x\cap U\neq \emptyset$.
We will consider separately two cases: when $Q$ contains a rational point and when it doesn't.
We will construct some $\gamma \in \Gamma$ with $\gamma^n x\in U$ for sufficiently big $n$. In the first case
$\gamma$ will be parabolic, and in the second case hyperbolic.

{\em Case 1.} Assume that $Q$ contains a rational point: let $u\in V_\bbQ$ be a non-zero vector with $q(u)=0$.
In this case, intersecting $Q$ with lines passing through $[u]\in Q$ we see that rational points are dense in $Q$.
Therefore we may assume that $[u]\in U$.

Let $G\subset O(V,q)$ be the unipotent subgroup with Lie algebra $u\wdg u^\perp\subset \Lambda^2 V$, where we identify the Lie algebra of the orthogonal group $O(V,q)$
with $\Lambda^2 V$ via the form $q$. For any $v\in u^\perp$, $w\in V$ and $t\in \bbR$ the action
of $e^{t u\wdg v}\in G$ on $w$ is given by the following formula:
\begin{equation}\label{eqn_action}
e^{t u\wdg v}w = w + tq(v,w)u - tq(u,w)v - \frac{t^2}{2} q(u,w)q(v,v)u.
\end{equation}
Rescaling $u$, we may assume that $u\in V_\bbZ$. The assumption $d\ge 2$ implies that we may find
$v\in u^\perp\cap V_\bbZ$ with $q(v,v)\neq 0$. Replacing $v$ by its multiple, we may
assume that the operator $\gamma = e^{u\wdg v}$ preserves $V_\bbZ$ and, furthermore, that $\gamma\in \Gamma$.
Choose $w\in V$ so that $x = [w]$. Then $q(u,w)\neq 0$ unless $u$ and $w$ are proportional.
In any case, we see from the equation (\ref{eqn_action}) that $\gamma^n x\to [u]$ for $n\to +\infty$, therefore $\gamma^n x \in U$ for $n$
big enough.

{\em Case 2.} Assume that $Q$ does not contain rational points. We may still find a two-dimensional rational
subspace $W_\bbQ\subset V_\bbQ$ such that $q|_{W_\bbQ}$ is non-degenerate, $x\notin \bbP(W)$ and $\bbP(W) \cap U\neq \emptyset$,
where $W = W_\bbQ\otimes \bbR\subset V$.
Then $Q\cap \bbP(W)$ consists of two points $[u]$ and $[u']$, where $u, u'\in W$ and we
may assume that $[u]\in U$.

By our assumption, the quadratic form $q|_W$ does not admit non-zero rational isotropic vectors.
After possibly rescaling $q$ by a rational factor, we may represent it in an appropriate basis by a diagonal matrix with entries $1$ and $-m$,
where $m>1$ is a square-free integer. We may identify the field $K = \bbQ[\sqrt{m}]$ with $W$
as a rational vector space, then $q$ is identified with the norm of an algebraic number.
It is well known that the group of units in the ring of algebraic integers $\OO_K$ is infinite.
Choose an element $a\in \OO_K^*$ of infinite order. Multiplication by $a$ preserves the
lattice $\OO_K\subset W$. Replacing $a$ by its power we may assume that there exists
an element $\gamma\in \Gamma$ such that $\gamma|_{W^\perp} = \mathrm{Id}$ and $\gamma|_W$ acts
as multiplication by $a$. Then $u$ and $u'$ are eigenvectors of $\gamma|_W$ with eigenvalues
$\lambda$ and $\lambda^{-1}$, and possibly replacing $\gamma$ by its inverse we may assume that $\lambda > 1$.
Recall that $x\notin \bbP(W)$ by construction and note that $x\notin \bbP(W^\perp)$,
since $q|_{W^\perp}$ is negative-definite. Since $\lambda$ is the
unique eigenvalue of $\gamma$ with absolute value greater than one, it is clear that
$\gamma^n x\to [u]$ as $n\to +\infty$ and so $\gamma^n x \in U$ for $n$ big enough.
\end{proof}

\hfill

Recall our notation: $Q = \{x\in \bbP(V)\st q(x) = 0\}$ is the boundary of $\bbP(V^+)$
and $\Gamma\subset O(V_\bbZ,q)$ is a subgroup of finite index.

\hfill

\proposition\label{prop_closure}
Assume that $d\ge 2$. For any $x\in \bbP (V)\setminus Q$ we have
$Q \subset\overline{\Gamma x}$, where the overline denotes the closure in $\bbP (V)$ in the analytic topology.

\hfill

\begin{proof}
The closed subset $\overline{\Gamma x} \subset \bbP(V)$ is $\Gamma$-invariant, so by \ref{prop_dense} it
is enough to show that $\overline{\Gamma x}$ contains at least one point of $Q$.
Choose $y\in V$ with $[y]=x$. We use the elements of $\Gamma$ constructed in the proof of \ref{prop_dense}.

{\em Case 1:} Assume that $Q$ contains a rational point. Then rational points are dense in $Q$, while
$Q\cap\bbP(y^\perp)$ is a proper closed subset of $Q$. Hence there exists $u\in V_\bbZ$ with
$q(u)=0$ and $q(u,y)\neq 0$. As in Case 1 of the proof of \ref{prop_dense}, choose
$v\in u^\perp\cap V_\bbZ$ with $q(v,v)\neq 0$ such that $\gamma=e^{u\wdg v}\in\Gamma$. By (\ref{eqn_action}),
$$
\gamma^m y = y + m\,q(v,y)u - m\,q(u,y)v - \frac{m^2}{2}\, q(u,y)q(v,v)\,u .
$$
Since $q(u,y)q(v,v)\neq 0$, the last term dominates as $m\to\infty$, so $\gamma^m x\to [u]\in Q$.

{\em Case 2:} Assume that $Q$ contains no rational points.
The two-dimensional subspaces of $V$ of signature
$(1,1)$ form a non-empty open subset of the Grassmannian. Rational subspaces are dense in it, and those
contained in $y^\perp$ form a proper closed subset. So we can choose a rational plane
$W$ of signature $(1,1)$ with $W\not\subset y^\perp$. Let
$\gamma\in\Gamma$ be as in Case 2 of the proof of \ref{prop_dense}: $\gamma|_{W^\perp}=\mathrm{Id}$, and
$\gamma|_W$ has isotropic eigenvectors $u$, $u'$ with eigenvalues $\lambda>1$ and $\lambda^{-1}$. Write
$y = \alpha u + \beta u' + z$ with $z\in W^\perp$. Then $(\alpha,\beta)\neq (0,0)$ and
$$
\gamma^m y = \lambda^m\alpha u + \lambda^{-m}\beta u' + z.
$$
If $\alpha\neq 0$, then $\gamma^m x\to[u]$ as $m\to+\infty$. Otherwise $\gamma^{m} x\to[u']$ as $m\to-\infty$.
In both cases $\overline{\Gamma x}\cap Q\neq\emptyset$.
\end{proof}

\hfill

\theorem\label{thm_boundary}
In the above setting let $\CC$ be one of the $\cS$-chambers for some $\cS$
satisfying (\ref{assumptions_S}) and $\overline{\CC}\subset V$ be the closure
of $\CC$ in $V$. Then the set $\bbP(\overline{\CC})\cap Q$ is nowhere dense in $Q$.

\hfill

\begin{proof}
Let $U\subset Q$ be an open subset. The set of positive vectors $V^+$ consists of two connected components and $\CC$ lies
in one of them. We can therefore find $u\in V$ such that its image $[u] \in \bbP(V)$ lies in $U$ and $q(w,u)>0$ for any $w\in\CC$.

By \ref{prop_closure}
for any $v\in \cS$ there exists a sequence $\gamma_i\in \Gamma$ such that $[\gamma_i v]\to [u]$ for $i\to +\infty$.
We may choose a sequence of vectors $v_i\in V^-$ such that $v_i$ is proportional to $\gamma_i v$ and
$v_i\to u$ for $i\to +\infty$. Let $H_i = v_i^\perp$, $H_i^- = \{w\in V\st q(w,v_i) < 0\}$
and $H_i^+ = \{w\in V\st q(w,v_i) > 0\}$. Since the vectors $v_i$ converge to $u$, the
hyperplanes $H_i$ converge to $u^\perp$. Since $\bbP(u^\perp)$ is tangent to $Q$ at
the point $[u]\in U$, it follows that infinitely many of the hyperplanes
$\bbP(H_i)$ have non-empty intersection with $U$. Therefore, passing to a subsequence,
we may assume that $\bbP(H_i^-)\cap U\neq \emptyset$ for all $i$.

Observe that for any $w\in \cap_i H_i^-$ we have $q(w,v_i) < 0$, so $q(w,u)\le 0$, because $v_i$ converge to $u$.
It follows that there must exist some $i_0$ such that $\CC\nsubseteq H_{i_0}^-$: otherwise for $w\in \CC$
we would have $q(w,u)\le 0$ which contradicts the choice of $u$.
Since $\cS$ is $\Gamma$-invariant, $\gamma_i v\in \cS$ for every $i$, so the chamber $\CC$ is either
contained in $H_i^-$ or in $H_i^+$. Therefore $\overline{\CC}\cap H_{i_0}^- = \emptyset$ and
we conclude that the non-empty open subset $U\cap \bbP(H_{i_0}^-) \subset U$ does not intersect $\bbP(\overline{\CC})$.
\end{proof}

\hfill

Recall that by our assumptions (\ref{assumptions_S}) we have $-\cS = \cS$.
Define
\begin{equation}\label{eqn_sprime}
\cS' = \{v\in \cS\st q(v,w)>0\,\, \forall w\in \CC\}
\end{equation}
and note that $\cS = \cS'\coprod (-\cS')$. For $v\in \cS'$
let $H_v^+ = \{w\in V\st q(v,w)>0\}$.

\hfill

\corollary\label{cor_closure}
In the above setting we have $\overline{\CC} = \overline{\cap_{v\in \cS'} H_v^+} = \cap_{v\in \cS'} \overline{H_v^+}$.

\hfill

\begin{proof}
Denote $\CC' = \cap_{v\in \cS'} H_v^+$.
Note that the set $\CC'\cap V^+$ is convex and
open, because the collection of hyperplanes $H_v$, $v\in \cS'$ is locally finite in $V^+$
by  \cite[Proposition 5.14]{SSV}.
It follows that $\CC = \CC'\cap V^+$, because $\CC$ is a connected component of the complement
to this collection of hyperplanes, and clearly $\CC \subset \CC'\cap V^+$.

Let us show that $\CC'\cap V^- = \emptyset$.
Assume that $v\in \CC'\cap V^-$. Since $\CC'$ is convex, it contains the segments connecting $v$
with the points of $\CC$, and the intersections of these segments with the boundary
of $V^+$ form a subset of the boundary of $\CC$ with non-empty interior.
This implies that the boundary of $\CC$ contains an open subset of the boundary of $V^+$,
contradicting \ref{thm_boundary}. It follows that $\CC'\cap V^- = \emptyset$.

We see that $\CC'\subset\overline{V^+}$. The convexity of $\CC'$ now implies that $\CC'\subset\overline{\CC}$.
Since $\CC\subset \CC'$, we also have $\overline{\CC}\subset\overline{\CC'}$ and the first
equality in the statement of the corollary follows.

To prove the second equality, note that any closed convex cone with non-empty interior
is the closure of its interior. Hence it suffices to prove that the cones
$\overline{\CC'}$ and $\cap_{v\in \cS'} \overline{H_v^+}$ have the same non-empty interior.
It is clear that the interior of $\overline{\CC'}$ is $\CC$, which is non-empty and is contained
in $\cap_{v\in \cS'} \overline{H_v^+}$. To prove the inclusion in the opposite direction,
note that a point in the interior of $\cap_{v\in \cS'} \overline{H_v^+}$ has an
open neighbourhood that lies in $H_v^+$ for all $v\in \cS'$, therefore lies in the
interior of $\CC'$. This completes the proof. 
\end{proof}

\hfill

Applying the above statements to the ample cone of a projective hyperk\"ahler manifold we recover the following theorem proved by
Kov\'acs in the K3 case, \cite{Ko1}, \cite{Ko2} (see also \cite{Mat} and \cite{_Denisi_}). Let $M$ be a hyperk\"ahler manifold.
Denote $V_\bbZ = \NS(M)$ and let $\Gamma\subset O(V_\bbZ, q)$ be the image of the Hodge monodromy group
under the natural representation, where $q$ is the BBF form.
By \cite[Lemma 6.23]{_Markman:survey_} $\Gamma$ is of finite index in $O(V_\bbZ, q)$.
Recall the definitions of the cones in the cohomology of $M$ given in Subsection \ref{sec_cones}.

\hfill

\theorem\label{thm_am}
Assume that $M$ is a projective hyperk\"ahler manifold with Picard number $\rho(M)>2$.
\begin{enumerate}
\item If $\NS(M)$ does not contain MBM classes, then $\AA_M = \CC_M\cap \NS_\bbR(M)$.
Otherwise the boundary of $\AA_M$ is nowhere dense in the boundary of $\CC_M\cap \NS_\bbR(M)$.
\item If $M$ does not contain any prime exceptional divisors then $\overline{\MM_M} = \overline{\CC_M}\cap \NS_\bbR(M)$.
Otherwise the boundary of $\overline{\MM_M}$ is nowhere dense in the boundary of $\CC_M\cap \NS_\bbR(M)$.
\end{enumerate}

\hfill

\begin{proof}
{\em Part 1.} Let $\cS = \MBM(M)\cap \NS(M)$. Since $\AA_M = \KK_M\cap \NS_\bbR(M)$,
we deduce the following from the description of $\KK_M$ given in Subsection \ref{sec_cones}: $\AA_M$ is a connected component of
the set $(\CC_M\cap\NS_\bbR(M)) \setminus \cup_{v\in \cS} v^\perp$. Since $\cS$ is $\Gamma$-invariant,
we may apply \ref{thm_boundary}.

{\em Part 2.} Analogously to part 1, but with $\cS$ defined to be the set of divisorial
MBM classes.
Theorem 6.17 in \cite{_Markman:survey_} shows that the interior of $\MM_M$
is one of the chambers in the decomposition of $\CC_M\cap\NS_\bbR(M)$ defined by $\cS$. We again
conclude by \ref{thm_boundary}.
\end{proof}

\hfill

We recall that $\mathrm{PE}(M)\subset \NS(M)$ denotes the set of classes of prime exceptional divisors.

\hfill

\theorem\label{thm_p}
Assume that $M$ is a projective hyperk\"ahler manifold with $\rho(M)>2$.
If $\mathrm{PE}(M)$ is empty, then $\PP_M = \overline{\CC_M}$, otherwise 
$$
\PP_M\cap \NS_\bbR(M) = \overline{\sum_{v\in \mathrm{PE}(M)} \bbR_{+}\cdot v}.
$$

\hfill

\begin{proof}
We use the description of $\PP_M$ given in (\ref{eqn_pseudoeffective}).
The case when $\PE(M)=\emptyset$ is clear, so we may assume the contrary.

Denote by $\cS$ the set of divisorial MBM classes and for $v\in \cS$
denote $H_v^+ = \{u\in \NS_\bbR(M)\st q(u,v)>0\}$. Then $\cS$ defines
the wall and chamber decomposition of $\CC_M\cap \NS_\bbR(M)$ and
we denote by $\CC$ the chamber whose closure is $\overline{\MM_M}$.
By \ref{cor_closure} we have $\overline{\CC} = \cap_{v\in \cS'}\overline{H_v^+}$,
where $\cS'$ is defined in (\ref{eqn_sprime}). On the other hand,
$\overline{\CC}$ is dual to $\PP_M\cap \NS_\bbR(M)$ and so we have
$\PE(M)\subset \cS'$, hence $\overline{\CC}\subset \cap_{v\in \PE(M)}\overline{H_v^+}$.
The latter inclusion must be equality, because the set of faces of
the cone $\CC$ inside the positive cone is given by the orthogonal
complements to the negative extremal rays of $\PP_M\cap \NS_\bbR(M)$,
i.e. by the vectors $v\in\PE(M)$.

We have shown that $\overline{\MM_M} =  \overline{\CC} = \cap_{v\in \PE(M)}\overline{H_v^+}$.
The claim now follows from duality between $\overline{\MM_M}$ and $\PP_M\cap \NS_\bbR(M)$.
\end{proof}

\subsubsection{The case when $\cS$ is contained in a proper subspace of $V$}

We go back to the setting introduced in the beginning of Subsection \ref{sec_hyperbolic}.
We additionally assume that there exists a proper subspace $N\subset V$ such that $\cS\subset N$.
Then $Q^\circ = Q\setminus \bbP(N)$ is an open subset of $Q$.
For $v\in \cS$ let $Q_v= Q\cap \bbP(H_v)$, where as above
$H_v = v^\perp$.

\hfill

\lemma\label{lem_finite_boundary}
The collection of subvarieties $Q_v$, $v\in\cS$ is locally finite in $Q^\circ$:
for any $x\in Q^\circ$ there exists an open subset $U\subset \bbP(V)$
such that $x\in U$ and $U$ intersects $Q_v$ only for finitely many $v\in \cS$.

\hfill

\begin{proof}
Let $u\in V$ be such that $[u]=x$. The idea is to translate the vector $u$
in the direction given by some $w\in N^\perp$, so that the translated
vector is contained in the positive cone. Local finiteness of the hyperplanes
$H_v$, $v\in \cS$ inside the positive cone is known, so we may find a ball
around the translated vector
that intersects only finitely many of the hyperplanes. Translating this ball back,
we get a neighbourhood of $u$ that we need. Note that $N^\perp\subset H_v$
for all $v\in \cS$, so that all hyperplanes are invariant under translation
in the direction of any $w\in N^\perp$.

More precisely, we note that since $x\in Q^\circ$,
we have $u\notin N$. The latter implies that there exists $w\in N^\perp$
with $q(u,w)\neq 0$. It follows that $u' = u+t_0w\in V^+$ for some $t_0\in \bbR$.
By \cite[Proposition 5.14]{SSV} the collection of hyperplanes $H_v$, $v\in\cS$
is locally finite in $V^+$, so there exists an open ball $B'\subset V^+$ containing $u'$
that intersects only finitely many hyperplanes $H_v$. Set $B = B'-t_0w$ and
note that since $w\in N^\perp$ the ball $B$ also intersects only finitely
many hyperplanes $H_v$, $v\in \cS$. Defining $U$ to be the image of $B$ in $\bbP(V)$
we conclude the proof of the claim.
\end{proof}

\hfill

Assume that $\CC\subset V^+$ is one of the $\cS$-chambers and 
let $\bnd = \bbP(\overline{\CC})\cap Q$ be the isotropic part of the boundary
of $\bbP(\CC)$. Denote by $\bndint$ the interior of $\bnd$ as a subset
of $Q$.

\hfill

\proposition\label{prop_chamber}
In the above setting the set $\bnd$ has the following properties.
\begin{enumerate}
\item[(1)] The intersection $\bndint\cap Q^\circ$ is non-empty and dense in $\bnd\cap Q^\circ$;
\item[(2)] If there exists $w\in N^\perp$ with $q(w)<0$ then $\bndint$ is dense in $\bnd$;
\item[(3)] If there exists a subspace $W\subset V$ containing $N$,
such that $\dim(W^\perp)\ge 2$ and $q|_{W^\perp}$ is negative definite,
then $\bndint$ and $\bnd$ are arcwise connected.
\end{enumerate}

\hfill

\begin{proof} {\em Proof of (1).} Note that when $q|_N$ is negative definite, then $\cS$ is finite, because it is
$q$-bounded and discrete. In this case the collection of hyperplanes $H_v$, $v\in \cS$
is finite and this easily implies part (1). It remains to prove (1) when $q|_N$ is either
degenerate or has mixed signature. Note that in both cases $q|_{N^\perp}$ is negative semi-definite.

Since $\CC\subset V$ is open and $N\subset V$ is a proper subspace,
there exists $u\in \CC\setminus N$. Since $u\notin N$, there exists $w \in N^\perp$
with $q(u,w)\neq 0$. By the remark above, we may assume that
$q(w)\le 0$. It follows that there exists $t_0\in \bbR$ such that $u'=u+t_0w$ is isotropic.
Since for any $v\in \cS$ and $t\in\bbR$ we have $q(u+tw,v)=q(u,v)\neq 0$, the affine line
$u + \bbR w$ does not intersect any of the hyperplanes $H_v$, $v\in\cS$. 
It follows that $\overline{\CC}$ contains the segment joining $u$ and $u'$, so $[u']\in \bnd\cap Q^\circ$.
Since $u$ lies in the interior of $\overline{\CC}$, we have $[u']\in \bndint\cap Q^\circ$,
so $\bndint\cap Q^\circ$ is non-empty. The fact that it is dense in $\bnd\cap Q^\circ$
follows from \ref{lem_finite_boundary}: any point $x\in \bnd\cap Q^\circ$ has an open
neighbourhood $U\subset \bbP(V)$ that intersects only finitely many hyperplanes $H_v$, $v\in\cS$.
This finite collection of hyperplanes and the quadric $Q$ cut $U$ into finitely many
connected components, and $\bbP(\CC)\cap U$ is one of these components. The claim of part (1)
follows easily from this.

{\em Proof of (2).} Let us fix a subspace $W$ such that $q|_{W^\perp}$
is negative definite. By assumption we may find such a $W$ with $W^\perp$ at least
one-dimensional. Then $N\subset W$ and therefore $\cS\subset W$.

Let $W^+ = W\cap V^+$ and $\CC_W = \CC\cap W.$ Note that the
linear projection of $V$ onto $W$ along $W^\perp$ induces a fibration $p\colon \bbP(\overline{V^+}) \to \bbP(\overline{W^+})$.
The fibre $F_p = p^{-1}[v]$ over a point $[v]\in \bbP(\overline{W^+})$ is isomorphic to the projectivization
of the cone
\begin{equation}\label{eqn_fibre}
\{av+w\st a\in \bbR,\, w\in W^\perp,\, q(av+w)\ge 0\}\subset W^\perp+\bbR v.
\end{equation}
This shows that $F_p$ is either a closed ball of dimension $\dim(V) - \dim(W)$
or a single point.
Since $\cS\subset W,$ each fibre $F_p$ is either contained in $\bbP(H_v)$ or does
not intersect it, for any given $v\in\cS$. It follows that $p$ induces a fibration
$\pi\colon \bbP(\overline{\CC}) \to \bbP(\overline{\CC_W})$ with the same fibres
as $p$. Note the following property of the fibration $p$: for any $x\in \bbP(W)\cap Q$
it is clear from (\ref{eqn_fibre}) that $p^{-1}(x) = \{x\}$. Moreover, if $y_i\in \bbP(\overline{V^+})$ is an arbitrary sequence
of points such that the sequence $\{p(y_i)\}$ converges to $x$, then the sequence $\{y_i\}$ also converges to $x$.
Analogous statements also hold true for the fibration $\pi$. 

For a point $x\in \bnd\setminus \bbP(N)$, we see from part (1) that $x$ lies
in the closure of $\bndint$.
To conclude the proof of part (2), take a point $x\in \bnd\cap\bbP(N)$ and a sequence of points
$y_i\in \bbP(\CC)$ converging to $x$. Since $x\in \bbP(W)$, we have $\pi(x)=x$
and we may choose arbitrary $z_i\in \bndint$ such that $\pi(z_i) = \pi(y_i) \in \bbP(\CC)$.
Then $\{z_i\}$ converges to $x$ by the remark at the end of the previous paragraph.

{\em Proof of (3).} We may now additionally assume that $\dim(W^\perp)\ge 2$. We use the
fibration $\pi$ constructed above. Let $x_1,x_2\in \bndint$. Then
$y_1 =\pi(x_1)$ and $y_2 = \pi(x_2)$ lie in $\bbP(W^+)$. Connect $y_1$ to $y_2$
by a segment and lift this segment to $\bndint$ by taking a section of $\pi$.
This gives an arc connecting $x_1$ to some point $x_2'\in \bndint$, where
$\pi(x_2')=\pi(x_2)$. The condition on the dimension of $W^\perp$ implies
that the boundary of the fibre $\pi^{-1}(y_2)$ is a sphere of positive dimension,
in particular it is connected. Therefore we may connect $x_2$ and $x_2'$ by an
arc lying in that sphere. This proves the claim for $\bndint$, and the
proof for $\bnd$ is analogous. 
\end{proof}

\hfill

We now apply the above proposition to the K\"ahler cone of a hyperk\"ahler manifold.

\hfill

\theorem\label{thm_kahler_boundary}
Let $M$ be a hyperk\"ahler manifold, $\CC_M\subset H^{1,1}(M, \bbR)$ its positive cone and $\KK_M$ its K\"ahler cone.
Let $\bnd_M = (\overline{\KK_M}\cap \partial\CC_M)\setminus\{0\}$ be the isotropic part of the boundary of $\KK_M$
and $\bnd_M^\circ = \bnd_M\setminus \NS_\bbR(M)$. Denote by $\rho(M)$ the Picard number of $M$.

\begin{enumerate}
\item If $\rho(M)\le h^{1,1}(M) - 1$ then the interior of $\bnd_M^\circ$ in $\partial\CC_M$ is non-empty and dense in $\bnd_M^\circ$;
\item If $\NS_\bbR(M)^\perp \subset H^{1,1}_\bbR(M)$ contains a $q$-negative vector
(for instance, if $\rho(M)\le h^{1,1}(M) - 2$) then the interior of $\bnd_M$ in $\partial\CC_M$ is non-empty and dense in $\bnd_M$;
\item If $\NS_\bbR(M)^\perp \subset H^{1,1}_\bbR(M)$ contains a $q$-negative subspace of dimension two
(for instance, if $\rho(M)\le h^{1,1}(M) - 3$) then both $\bnd_M$ and its interior in $\partial\CC_M$
are arcwise connected.
\end{enumerate} 

\hfill

\begin{proof}
This follows directly from \ref{prop_chamber} applied to $V = H^{1,1}_\bbR(M)$, $N = \NS_\bbR(M)$,
$\cS = \MBM(M)\cap\NS(M)$, $\CC = \KK_M$.
\end{proof}

\section{The monodromy group and automorphisms of hyperk\"ahler manifolds}
\label{_monodromy_Section_}

Here we present some basic results about the rational
curves, automorphisms and the K\"ahler cone of a
hyperk\"ahler manifold. We follow
\cite{_Verbitsky:Torelli_,_AV:MBM_,_AV:automorphisms_,_AV:hyperb_}.

\subsection{The monodromy orbifold of a hyperk\"ahler manifold}

We recall that an {\bf orbifold} is a
topological space with a maximal atlas consisting
of quotients of open subsets of $\bbR^n$ by finite groups,
with smooth transition functions respecting the group
actions, see \cite{MP} for a precise definition
and \cite{ALR} for a detailed discussion. It is sometimes
useful to think of orbifolds in terms of topological
groupoids, see loc. cit. for this point of view.
The basic example (and the only one that we will need)
is the {\bf action groupoid}: if $M$ is a smooth manifold
and $\Gamma$ is a discrete group acting on $M$ properly discontinuously,
then the orbifold groupoid is given by the two maps $s,t\colon \Gamma\times M\to M$,
where $s$ is the projection to $M$ and $t$ is the map defining the action,
see \cite[section 1]{ALR} for details. The orbifold represented by
the action groupoid is the quotient orbifold $M/\Gamma$.

A {\bf covering} of an orbifold is a 
map $\pi\colon X_1 \to X_2$ which can be locally
represented as $B/\Gamma_1 \arrow B/\Gamma_2$,
where $B\subset \bbR^n$ is an open ball and $\Gamma_1 \subset \Gamma_2$,
and such that any point in $X_2$ has a neighbourhood
$U_2$ such that its preimage is a union of
open subsets in $X_1$ mapping to $U_2$ in this
fashion. An orbifold $X$ is called {\bf simply connected}
if any connected covering of $X$ is isomorphic to $X$.
A {\bf universal covering} of $X$ is a connected
covering which is simply connected. One can show
that any orbifold admits a universal covering which is unique up to isomorphism.
The {\bf fundamental group} of a connected
orbifold $X$ is the group of automorphisms
of its universal covering compatible with
the projection to $X$. See \cite[section 2.2]{ALR}
for the definition of the orbifold fundamental group
in terms of groupoids. 

Going back to the example of the quotient orbifold,
let $\Gamma$ be a discrete group acting properly
discontinuously on a contractible manifold $Z$.
Then $Z$ is the universal covering of the
orbifold $Z/\Gamma$, and $\pi_1(Z/\Gamma)= \Gamma$.

Let $\Gamma \subset SO^+(1, n)$
be a lattice, i.e. a discrete subgroup of finite co\-volume.
A {\bf hyperbolic orbifold} is the quotient ${\Bbb H}^n/\Gamma$ of the
hyperbolic space.
It is called a {\bf hyperbolic manifold} if $\Gamma$ has
no elements of finite order (then, by the classification of isometries of the hyperbolic space, $\Gamma$ acts freely on $\Bbb H$). 

Let $G$ be a semisimple algebraic group defined over $\Q$.
An {\bf arithmetic lattice} is a subgroup $\Gamma \subset G$
which is commensurable with $G_\Z$. Recall that by
a theorem of Borel and Harish-Chandra \cite{_Morris:Ratner_},
for an algebraic Lie group $G$ defined over $\Q$
and without non-trivial rational characters 
(such as a semisimple group), $G_\Z$ has finite covolume.
This theorem implies
that for an arithmetic lattice $\Gamma\subset SO^+(1,n)$
the quotient ${\Bbb H}^n/\Gamma$ is a hyperbolic orbifold.

We will consider hyperbolic orbifolds that arise from the
action of the automorphism group of a hyperk\"ahler manifold
on its cohomology. For a hyperk\"ahler manifold $M$
we will denote by $\Diff_0(M)$ the connected component of
its diffeomorphism group containing the identity (the group of isotopies).
We will denote by $\Comp(M)$ the (infinite-dimensional) space
of all complex structures of hyperk\"ahler type on $M$, and we define 
the {\bf Teich\-m\"uller space} as $\Teich(M)=\Comp(M)/\Diff_0(M)$.
It is known that $\Teich(M)$ is a finite-dimensional non-Hausdorff
complex analytic space (see e.g. \cite{_Catanese:moduli_}, \cite{_Verbitsky:Torelli_}).

Let $\MCG:=\Diff(M)/\Diff_0(M)$ denote the mapping class group of $M$.
For a complex structure $I\in\Comp(M)$ let 
$\Gamma^I$ be the  subgroup of $\MCG$ preserving the
connected component $\Teich^I\subset \Teich$. The
image of $\Gamma^I$ in $O(H^2(M, \Z))$ is
called {\bf the monodromy group} of $(M,I)$,
denoted by $\Mon(M,I)$ or $\Mon(M)$ when $I$ is clear from the context.

\hfill

\remark
This notion is due to E. Markman (\cite{_Markman:survey_}).
The monodromy group can be alternatively defined
as the subgroup of $O(H^2(M, \Z))$ generated by
the monodromy operators for all Gauss--Manin local systems
obtained from deformations of $(M,I)$.

\hfill

Let $I\in\Comp(M)$ be a complex structure on $M$.
The {\bf Hodge mono\-dromy group} is the group $\Mon_I$
of all $s\in\Mon(M,I)\subset O(H^2(M, \Z))$ that preserve the Hodge
decomposition $H^2(M, \C)= H^{2,0}(M) \oplus H^{1,1}(M)\oplus H^{0,2}(M)$.

We recall (see \cite{_Markman:survey_}) that
the automorphism group of $(M,I)$ is discrete, the restriction map
$\Aut(M,I) \to \Aut(M,I)|_{H^2(M, \Z)}$
has finite kernel, and its image is the set of all
$\tau\in \Mon_I$ that preserve the K\"ahler cone. Slightly more generally, the following claim is true.

\hfill

\claim\label{chambers} 
The orthogonal hyperplanes to the MBM classes partition the positive cone $\CC_M$ into chambers, one of which is the K\"ahler cone
$\KK_M$. The group $\Mon_I$ acts on the set of the chambers (called {\bf K\"ahler chambers}), and its elements which lift to automorphisms of 
$M$ are exactly the ones which preserve $\KK_M$. The K\"ahler chambers are the K\"ahler cones of the birational models of $M$ and their monodromy images. 

\hfill

We sometimes denote by $H^{1,1}(M,\Z)$ 
the group $H^2(M, \Z)\cap H^{1,1}(M)=\NS(M)$.
Clearly, $\Mon_I$ acts on $H^{1,1}(M,\Z)$.

The following theorem is well-known; for the convenience of the reader, we highlight a part of the proof.

\hfill

\theorem\label{_Mon_I_arithm_Proposition_}
Let $M$ be a projective hyperk\"ahler manifold with complex
structure $I$. Then
the natural homomorphism $\Mon_I\to O(H^{1,1}(M,\Z))$
has finite kernel, and its image is of finite index
in $O(H^{1,1}(M,\Z))$.

\hfill

\pstep
As shown in \cite{_Verbitsky:Torelli_}, 
the group $\Mon$ has finite index in $O(H^2(M, \Z))$.
Therefore, $\Mon_I$ has finite index in the group $O_I(H^2(M, \Z))$
of all $s\in O(H^2(M, \Z))$ preserving the Hodge
decomposition. Let $\Lambda\subset H^2(M,\bbZ)$ be the
orthogonal complement of $H^{1,1}(M,\Z)$ in $H^2(M,\Z)$.
Since we assume that $M$ is projective, the lattice
$H^{1,1}(M,\Z)\oplus \Lambda$ has finite index in $H^2(M, \Z)$.
It follows that the image of $O_I(H^2(M, \Z))$
in $O(H^{1,1}(M,\Z))$ has finite index.
This implies that the image of 
$\Mon_I$ in $O(H^{1,1}(M,\Z))$
has finite index.

\hfill

{\bf  Step 2:} 
It remains to show that the kernel
of the natural map $\Mon_I\to O(H^{1,1}(M,\Z))$
is finite. This would follow if we prove that
the kernel $K$ of the natural map $O_I(H^2(M, \Z))\arrow O(H^{1,1}(M,\Z))$
is finite. Since $M$ is projective and $K$ acts trivially
on $H^{1,1}(M,\bbZ)$, the group $K$ is contained in the group
of automorphisms of the integral polarized Hodge structure $H^2(M,\bbZ)$.
The latter group is finite. 
\endproof

\hfill


The intersection
\[ \Pos_M= \CC_M\cap \big(H^{1,1}(M,
\Z)\otimes_\Z \R\big)
\]
will be called {\bf the rational positive cone}.
Recall that the intersection $\Pos_M\cap\, \KK_M$ is the
ample cone $\AA_M$. Clearly, we have an exact analogue of \ref{chambers} in this setting. The resulting chambers,
one of which is $\AA_M$, are called {\bf ample chambers}. 
The projectivization ${\Bbb P}\Pos_M$
is a hyperbolic space; assuming that $M$ is projective,
it follows from \ref{_Mon_I_arithm_Proposition_} that
the group $\Mon_I$ acts on ${\Bbb P}\Pos_M$ as an arithmetic lattice. 

\hfill

\definition\label{def_ample_poly}
The quotient orbifold ${\Bbb P}\Pos_M/\Mon_I$
will be called {\bf the mo\-nodromy orbifold} of the hyperk\"ahler 
manifold $(M,I)$. Define {\bf the hyperbolic ample polyhedron}
as the image of ${\Bbb P}\AA_M$ in the monodromy orbifold.

\hfill

\remark\label{selberg} The monodromy orbifold has a finite covering which is a manifold. Indeed by Selberg's lemma the image of $\Mon_I$ in $O(\NS(M))$ has a torsion-free subgroup of finite index, and elements with fixed points on ${\Bbb P}\Pos_M$ are torsion.
The same applies to any hyperbolic orbifold ${\Bbb H}/\Gamma$ (recall that $\Gamma$ is a lattice in a semisimple group $SO(1,n)$). For our purposes, such as studying limit sets, we may as well work with 
manifolds.




\subsection{Shape of the K\"ahler cone}

Consider the K\"ahler cone $\KK_M$ of a hyperk\"ahler manifold $M$.
The K\"ahler cone is open and convex in $H^{1,1}(M,\bbR)$.
A {\bf wall} of $\KK_M$ is the intersection
of the boundary $\6\KK_M$ and a hyperplane $W\subset H^{1,1}(M,\R)$
that has a non-empty interior (in both). A {\bf face} of $\KK_M$ is a non-empty intersection of walls; in particular, walls are faces of maximal dimension.
A hyperplane $W$ supporting a wall of $\KK_M$ is the orthogonal hyperplane to an MBM class. So a face of $\KK_M$ is obtained by intersecting the boundary of $\KK_M$ and a finite number of orthogonal hyperplanes to MBM classes.

The walls and faces are defined in the same way for the ample cone $\AA_M$ of a projective hyperk\"ahler manifold $M$, by considering 
$\NS_{\bbR}(M)$ instead of $H^{1,1}(M,\bbR)$. 

In
\cite{_AV:Kaw_Mor_}, it has been proved that the Hodge monodromy group acts with
finitely many orbits on the set of MBM classes of type
(1,1). As a consequence, the following version of the
Kawamata--Morrison conjecture has been obtained.

\hfill

\theorem \cite{_AV:hyperb_}\label{_faces_Kaw_Mor_Theorem_}
Let $(M,I)$ be a hyperk\"ahler manifold. Then 
the set $F$ of walls of the K\"ahler cone $\KK_M$ is
locally finite outside of the boundary of the positive cone, 
and the group of complex automorphisms
$\Aut(M,I)$ acts on $F$ with finitely many orbits. If $M$ is projective, the same holds for the ample cone.

\hfill

\remark 
It follows from the above theorem that the boundary $\6\KK_M$
of the K\"ahler cone is the union of the faces 
and ``the parabolic part'', i.e. the set
$\{\eta \in \6\KK_M \st q(\eta)=0\}$.

\hfill

\definition
Let $V$ be a quadratic vector space of signature $(1,n)$,
and ${\Bbb H}^n= {\Bbb P}V^+$ the corresponding hyperbolic space.
{\bf A hyperbolic hyperplane}
in the hyperbolic space ${\Bbb H}^n= {\Bbb P}V^+$ is
the image of ${\Bbb P}W^+$, where $W\subset V$
is a hyperplane of signature $(1, n-1)$.
Note that the set of hyperbolic hyperplanes
is in bijective correspondence with
${\Bbb P}V^-$, where $V^-$ is the 
set of negative vectors. {\bf A locally geodesic  hypersurface}
in a hyperbolic orbifold ${\Bbb H}^n/\Gamma$ is
the image of a hyperplane ${\Bbb P}W^+\subset {\Bbb H}^n$.
Note that a locally geodesic hypersurface, if it is closed,
is also a hyperbolic orbifold.

\hfill

\remark A locally geodesic hypersurface is a closed subset
of ${\Bbb H}^n/\Gamma$ if and only if $\Gamma\cap SO(1, n-1)$ is a lattice in $SO(1,n-1)$ (the stabilizer of $W$).
This is a consequence of Ratner theory when $n>2$ (see \ref{_Ratner_closure_Remark_}), and elementary when $n=2$. In the presence of an integral structure, that is, when 
${\Bbb H}^n={\Bbb P} V^+$ with $V=V_\Z\otimes \R$, and $\Gamma$ is arithmetic, this amounts to saying that $W$ is rational. Of course the same applies to the case when $W\subset V$ is a subspace of higher codimension $k$ (but still of signature $(1, n-k)$). One calls ${\Bbb P}W^+$ {\bf a hyperbolic subspace}, and its image in 
${\Bbb H}^n/\Gamma$ {\bf a locally geodesic (immersed) submanifold}.

\hfill

\definition\label{_polyhe_orbifo_Definition_}
Let ${\cal H}= {\Bbb H}^n/\Gamma$
be a hyperbolic orbifold, and $S_1, \ldots, S_m\subset {\cal H}$
be a collection of closed locally geodesic hypersurfaces.
A {\bf convex hyperbolic orbifold with 
 polyhedral boundary} (CHOPB)
is a connected component of the complement
${\cal H}\backslash \bigcup_i S_i$.

\hfill

\theorem\label{_kahler_convex_hy_polyhe_Theorem_}
Let $(M,I)$ be a projective hyperk\"ahler manifold,
and ${\cal P}:= {\Bbb P}\AA_M/\Aut(M,I)$
be its hyperbolic ample polyhedron.
Then ${\cal P}$ is a convex hyperbolic orbifold with 
polyhedral boundary,
and 
$\pi_1({\cal P})$ is equal to the image of
$\Aut(M,I)$ in $O(H^{1,1}(M,\Z))$. Moreover,
the natural surjection $\Aut(M,I) \to \pi_1({\cal P})$
has finite kernel.

\hfill

\proof 
The orthogonal hyperplanes to the MBM classes of type $(1,1)$ partition ${\Bbb P}\Pos_M$ into chambers, one of which is ${\Bbb P}\AA_M$.
Let ${\cal M}= {\Bbb P}\Pos_M/\Mon_I$
be the monodromy orbifold of $M$.
The preimage of ${\cal P}\subset {\cal M}$
in the universal cover ${\Bbb P}\Pos_M$ of 
${\cal M}$ is the union of ${\Bbb P}\AA_M$ and some other chambers which are its $\Mon_I$-translates, and ${\Bbb P}\AA_M$ is clearly 
the universal covering of ${\cal P}$. 
Let $S_i\subset {\cal M}$ be the closed hypersurfaces obtained
as images of orthogonal complements to the
MBM classes of type $(1,1)$. Then ${\cal P}$ is a connected component
of the set ${\cal M}\backslash \bigcup S_i$.
 As follows from the
solution of the Kawamata--Morrison conjecture
(\cite{_AV:Kaw_Mor_}), the set $\{S_i\}$ 
is finite. Finally, 
$\pi_1({\cal P})$ is the subgroup of $\Mon_I = \pi_1(\bbP\Pos_M/\Mon_I)$
which preserves the ample cone $\AA_M \subset \Pos_M$.
By \cite{_Markman:survey_}, this is precisely the
image of $\Aut(M,I)$ in $O(H^{1,1}(M,\Z))$, and the kernel is finite by \ref{_Mon_I_arithm_Proposition_}.
\endproof


\section{Limit sets of convex hyperbolic orbifolds with polyhedral boundary}
\label{_limit_points_Section_}

By \ref{_kahler_convex_hy_polyhe_Theorem_}, the automorphism group of a projective
hyperk\"ahler manifold is, up to a finite kernel, the orbifold fundamental group
of a convex hyperbolic orbifold with polyhedral boundary (CHOPB). In this section
we study the limit sets of such groups in general.


\subsection{The absolute and the limit 
set of a group acting by hyperbolic isometries}

\definition
Let $V$ be a real vector space equipped with 
a scalar product $q$ of signature $(1, n)$,
and $V^+$ the set of all vectors with positive square.
Then ${\Bbb H}^n= {\Bbb P} V^+$ is the hyperbolic space
of constant negative curvature.
{\bf The absolute} $\Abs$ is 
the projectivization of the set of all vectors
with square 0; it is identified with the boundary
of ${\Bbb H}^n$. Clearly, $\Abs$ is diffeomorphic to 
the sphere $S^{n-1}$, and the union ${\Bbb H}^n\cup \Abs$ is compact.

%
%
%
%
%

\hfill

\claim\label{_geode_distance_Claim_}
\begin{enumerate}
\item[(i)] Let $a_i, b_i\in{\Bbb H}^n$ be sequences with $a_i\to x\in\Abs$ and $b_i\to y\in{\Bbb H}^n\cup\Abs$,
$y\neq x$. Then $d(a_i,b_i)\to\infty$.
\item[(ii)] Let $y_i\in{\Bbb H}^n$ converge to $x\in\Abs$, and let $p_i\in{\Bbb H}^n$ satisfy $d(p_i,y_i)\le R$
for some constant $R$. Then $p_i\to x$.
\end{enumerate}

\hfill

\proof
We use the hyperboloid model: ${\Bbb H}^n=\{a\in V\ |\ q(a)=1,\ q(a,e)>0\}$ for a fixed $e$ with $q(e)=1$,
where $\cosh d(a,b)=q(a,b)$. Fix the Euclidean norm $\|w\|^2 = 2q(w,e)^2-q(w)$. For $a\in{\Bbb H}^n$ we have
$q(a,e)\ge 1$, hence $\|a\|\ge 1$.

(i) Write $a_i=t_i\xi_i$ and $b_i=s_i\eta_i$ with $\|\xi_i\|=\|\eta_i\|=1$, so $s_i\ge 1$. Passing to a subsequence,
$\xi_i\to\xi$ and $\eta_i\to\eta$ with $[\xi]=x$ and $[\eta]=y$. Since $q(\xi_i)=t_i^{-2}\to q(\xi)=0$, we get $t_i\to\infty$.
Now $\xi$ and $\eta$ lie in the closure of the same component of $V^+$, $\xi$ is isotropic, and $\eta$ is not
proportional to $\xi$, hence $q(\xi,\eta)>0$. Therefore
$\cosh d(a_i,b_i) = t_is_i\,q(\xi_i,\eta_i)\ge t_i\,q(\xi_i,\eta_i)\to\infty$.
Since this holds for a subsequence of any sequence, $d(a_i,b_i)\to\infty$.

(ii) Otherwise, by compactness of ${\Bbb H}^n\cup\Abs$, a subsequence of $p_i$ converges to some $y\neq x$.
Then $d(p_i,y_i)\to\infty$ along this subsequence by (i), a contradiction.
%
\endproof

%

\hfill

\claim
Let $\Gamma\subset SO^+(1,n)$
be a group of hyperbolic isometries
acting on ${\Bbb H}^n$, and
$a, b\in {\Bbb H}^n$
any two points. Let $\Lambda_{a}, \Lambda_b \subset \Abs$ 
be the set of all points in $\Abs$
obtained as accumulation points
for $\Gamma\cdot a, \Gamma\cdot b$.
Then $\Lambda_{a}=\Lambda_{b}$.

\hfill

\proof
Let $x\in \Lambda_{a}$, and $\{y_i=\gamma_i a \}\in \Gamma\cdot a$
be a sequence converging to $x$. Since $\gamma_i$ are isometries,
$d(\gamma_i b, \gamma_i a) =d(a,b)$, and
the sequence $\gamma_i b$ converges to $x$ as well.
Indeed, since ${\Bbb H}^n\cup \Abs$ is compact,
otherwise we would have a point $y\neq x$ which is
a limit of a subsequence $\gamma_i b$.
This is impossible by \ref{_geode_distance_Claim_} (i),
because $d(\gamma_i a,\gamma_i b)=d(a,b)$ is constant.
\endproof

\hfill

\definition
Under these assumptions,
{\bf the limit set} of the action of $\Gamma$ on ${\Bbb H}$ is
$\Lambda_{a}$, for some (and, therefore, any) $a\in {\Bbb H}$.

\hfill

\remark
Let $\Lambda$ be the limit set of $\Gamma$. Clearly,
$\Gamma$ acts on $\Abs\backslash \Lambda$ properly 
discontinuously.

\subsection{Convex hyperbolic orbifolds with 
polyhedral boundary in arithmetic hyperbolic orbifolds}





\definition
Let $V_\Z$ be a lattice
equipped with a scalar product of signature 
$(1,n)$, $V=V_\Z\otimes \R$, and $\Gamma \subset SO^+(V)$
a subgroup commensurable with $SO(V_\Z)$.
We call ${\Bbb P}V^+/\Gamma$, that is, the quotient of the hyperbolic space
by an arithmetic lattice $\Gamma$, {\bf an arithmetic hyperbolic orbifold}.
We sometimes use the notation $\Pos(V)$ for $V^+$.

\hfill

\remark\label{_Ratner_closure_Remark_}
Let $\Gamma\subset SO^+(1,n)$ be a discrete group of isometries
of finite covolume (that is, a lattice), 
and ${\Bbb P}W^+\subset {\Bbb P}V^+$
a hyperbolic subspace. Clearly,
the image of ${\Bbb P}W^+$ in
${\Bbb P}V^+/\Gamma$ has the same volume
as ${\Bbb P}W^+/\Gamma_W$, where $\Gamma_W$
is $\{\gamma\in \Gamma \ \ | \ \ \gamma(W)=W\}$.
Any integer lattice in $O(1,n)$ has finite covolume.
By Ratner's orbit closure theorem \cite{_Ratner:Raghunathan_},
in the form proved for hyperbolic manifolds by Shah \cite{_Shah:closures_},
the closure of the image of a hyperbolic subspace ${\Bbb P}W^+$ with $\dim W\geq 3$
in ${\Bbb P}V^+/\Gamma$ is the image of ${\Bbb P}U^+$ for a subspace
$U\supset W$ such that $\Gamma_U$ is a lattice in $O(U)$.
When $\Gamma$ is commensurable with $SO(V_\Z)$, the group $\Gamma_U$
is a lattice in $O(U)$ if and only if $U$ is rational (by the Borel
density theorem and the Borel--Harish-Chandra theorem).
In particular, any closed totally geodesic immersed submanifold of dimension at least 2
in ${\Bbb P}V^+/\Gamma$ is equal to ${\Bbb P}W^+/\Gamma_W$,
for some rational subspace $W\subset V $.
This implies the following description of
a convex hyperbolic orbifold with 
polyhedral boundary (CHOPB) in arithmetic
hyperbolic orbifolds.

\hfill

\claim
Let $V_\Z$ be a lattice
equipped with a scalar product of signature 
$(1,n)$, and $\Gamma \subset SO^+(V)$
a subgroup commensurable with $SO(V_\Z)$.
Consider a finite collection of
rational hyperplanes $W_i\subset V$ of 
signature $(1,n-1)$, let ${\goth S}$
be $\bigcup_i \Gamma W_i$, and 
$P$ a connected component
of ${\Bbb P} V^+\backslash {\Bbb P} {\goth S}$.
Denote by $\Gamma_P$ the group  $\{\gamma\in \Gamma \ \ | \ \ \gamma(P)=P\}$.
Then $P/\Gamma_P$ is a convex hyperbolic orbifold with 
polyhedral boundary, and all convex hyperbolic orbifolds with 
polyhedral boundary in ${\Bbb P} V^+/\Gamma$ are obtained this way.
\endproof

\hfill

The following question is crucial for this paper.

\hfill

\question
 What is the limit set of $\Gamma_P$ acting on $P$?

\hfill

The reason is the following theorem, true for all
hyperbolic orbifolds, but it is much easier
to state and prove it when  ${\cal H}$ is
compact. In Subsection \ref{_abs_of_polyhe_Subsection_}
we will explain how to generalize it
to all hyperbolic ${\cal H}$.

\hfill

\theorem\label{_limit_points_H_compact_Theorem_}
Let ${\cal H}:={\Bbb P} V^+/\Gamma$
be a hyperbolic orbifold, and ${\cal P}\subset {\cal H}$ 
a convex hyperbolic orbifold with 
polyhedral boundary. Let $P$ be a connected
component of the preimage of ${\cal P}$
in ${\Bbb P} V^+={\Bbb H}$; clearly,
$P$ is a convex polyhedron in ${\Bbb H}$
with hyperbolic faces. Denote by $\Abs P$
the set $\Abs\cap \bar P$.
Assume that ${\cal H}$ is compact.
Then $\Abs P$ is the limit set of 
$\Gamma_P$ acting on $P$.

\hfill

\pstep
Clearly, for any $x\in P$, its $\Gamma_P$-orbit
belongs to $P$, hence its limit set
belongs to $\Abs P$.

\hfill

{\bf Step 2:}
Conversely, 
let $x\in \Abs P$,
and let $l$ be a geodesic 
with an end in $x$. We may choose $l$ that contains a point in $P$.
Since $P$ is convex, a ray $l_+\subset l$ converging to $x$
also belongs to $P$. Choose a fundamental domain
$D$ of $\Gamma_P$-action on $P$.
Since $P/\Gamma_P$ is a closed subset of ${\cal H}$, and the latter
is assumed compact, $D$ can be chosen relatively compact.
Let $R:= \diam D$.

\hfill

{\bf Step 3:}
Choose $z\in D$. Every point $y\in P$ lies in $\gamma\bar D$ for some $\gamma\in\Gamma_P$,
hence $d(y,\gamma z)\le R$. Choose a sequence of points $y_i\in l_+$ converging to $x\in\Abs P$,
and $\gamma_i\in\Gamma_P$ with $d(\gamma_i z, y_i)\le R$.
By \ref{_geode_distance_Claim_} (ii), $\gamma_i z\to x$,
hence $x$ belongs to the limit set of $\Gamma_P$.
\endproof

\subsection{The thick-thin decomposition}

Here we give a brief introduction to cusp
points on hyperbolic orbifolds. For 
more details, see \cite{_Thurston:thick-thin_,_Kapovich:Kleinian_}.

\hfill

\definition\label{defn_horosphere}
{\bf A horosphere} in a hyperbolic space is 
a (punctured) sphere which is everywhere orthogonal to the pencil
of geodesics passing through a given point at infinity, and 
{\bf a horoball} is a ball bounded by a horosphere.
{\bf A cusp point} of an $n$-dimensional hyperbolic
orbifold ${\Bbb H}/\Gamma$ is a point $p\in\Abs{\Bbb H}$
whose stabilizer $\Gamma_p\subset \Gamma$ contains a free
abelian group of rank $n-1$ consisting of parabolic elements;
equivalently, $\Gamma_p$ is a cocompact lattice in the group
of isometries preserving the horospheres centred at $p$.
The stabilizer $\Gamma_p$ of a cusp point is called
{\bf a maximal parabolic} subgroup of $\Gamma$.
Horospheres centred at $p$ are the orbits of the unipotent
radical of the stabilizer of $p$ in $SO^+(1,n)$.

\hfill

\claim
Let $V$ be a quadratic vector space of signature
$(1,n)$, and $\Gamma\subset SO^+(V)$ a discrete
subgroup of finite covolume.
Consider a point $p\in \Abs$, and let 
$\Gamma_p\subset \Gamma$ be its stabilizer. Clearly,
$\Gamma_p$ acts on any horosphere $S$ tangent to the
boundary at $p$ by isometries.
In such a situation, $p$ is a cusp point 
if and only if $(S\backslash \{p\})/\Gamma_p$
is compact. \endproof

\hfill

\remark
A cusp point $p$ yields a {\bf cusp} 
in the quotient ${\Bbb H}/\Gamma$, that is,
a geometric end of ${\Bbb H}/\Gamma$ 
of the form $B/\Gamma_p$, where 
$B\subset \H$ is a horoball tangent
to the boundary at $p$ and $\Gamma_p$ contains ${\Bbb Z}^{n-1}$ as a subgroup of finite index.

\hfill

\theorem\label{_thick_thin_Theorem_}
 (Thick-thin decomposition) \\
Any $n$-dimensional complete hyperbolic 
orbifold $\H/\Gamma$ of finite volume can be represented as
a union of a ``thick part'', which is a compact orbifold
with boundary, and a ``thin part'', which is
a finite union of quotients isometric to $B/\Gamma_p$, where $B$ is a horoball
tangent to the boundary at a cusp point $p$, and
$\Gamma_p$ its stabilizer in $\Gamma$. Moreover,
$\Gamma_p$ contains $\Z^{n-1}$ as a subgroup of finite index, and
$B/\Gamma_p$ is diffeomorphic to a product of a finite quotient of a torus $T^{n-1}$
and an interval.

\hfill

\proof See 
\cite[Section 5.10]{_Thurston:thick-thin_} or
\cite[page 491]{_Kapovich:Kleinian_}. \endproof

\hfill

When the group $\Gamma\subset SO^+(V)$
is arithmetic, that is, commensurable
with a lattice $SO(V_{\Z})$,
the cusp points can be identified explicitly.
The following theorem seems to be well known.

\hfill

\theorem 
Let $V_\Z$ be a lattice
equipped with a scalar product of signature 
$(1,n)$, and $\Gamma \subset SO^+(V)$
a subgroup commensurable with $SO(V_\Z)$.
Then the cusp points of the hyperbolic
orbifold ${\Bbb H}/\Gamma$ are in 
bijective correspondence with $\Gamma$-orbits
on the set $\Abs \cap {\Bbb P} (V_\Q)$
of all rational points on $\Abs$.

\hfill

\proof From 
\cite[Proposition I.4.9]{_Borel_Ji_}
 (see also \ref{_thick_thin_Theorem_} or \cite{_Kapovich:Kleinian_})
we obtain that the cusps of a hyperbolic orbifold ${\cal H}$
correspond to maximal parabolic 
subgroups  $\Gamma_0\subset \pi_1({\cal H})=\Gamma$,
see \ref{defn_horosphere}.
The fixed point $\alpha$ of $\Gamma_0$-action on the
absolute is the corresponding cusp point of 
${\Bbb  H}/\Gamma$. This implies that
the line $\R \alpha$ coincides with $\cap_{\gamma\in \Gamma_0}\ker (\Id-\gamma)$,
which is rational, being the intersection of rational subspaces.

Conversely, consider a rational point $\alpha$ in the
absolute, and let $G_0\subset SO(1,n)$
be the maximal parabolic group fixing $\alpha$.
Since $G_0$ is a rational subgroup in $SO(1,n)$,
its intersection with the lattice $\Gamma$, commensurable
to $SO(V_\Z)$, is a lattice in $G_0$.
 \endproof

\subsection{The absolute of a universal covering of a CHOPB}
\label{_abs_of_polyhe_Subsection_}

Having defined the cusp points and introduced the
thick-thin decomposition, we can state and prove
the general version of \ref{_limit_points_H_compact_Theorem_}.

\hfill

\theorem\label{_Limit_set_Theorem_}
Let ${\cal H}:={\Bbb P} V^+/\Gamma$
be a hyperbolic orbifold, and ${\cal P}\subset {\cal H}$ 
a convex hyperbolic orbifold with 
polyhedral boundary. Let $P$ be a connected
component of the preimage of ${\cal P}$
in ${\Bbb P} V^+ ={\Bbb H}$; clearly,
$P$ is a convex polyhedron in ${\Bbb H}$
with hyperbolic faces. Denote by $\Abs P$
the set $\bar P\cap \Abs$.
Then $\Abs P$ is the union of the set of
all cusp points of ${\cal H}$ which belong to $\Abs P$ and the 
limit set of $\Gamma_P$ acting on $P$.

\hfill

\pstep
Clearly, for any $x\in P$, its $\Gamma_P$-orbit
belongs to $P$, hence its limit set
belongs to $\Abs P$.

\hfill

{\bf Step 2:} Consider the thick-thin decomposition
of ${\cal H}$, with the thick part ${\cal H}_0$,
and let $P_0$ be the intersection of $P$ and
the preimage of ${\cal H}_0$. Then 
$P_0$ is obtained by removing from $P$ a countable number
of horoballs; this set is $\Gamma_P$-invariant. 
A fundamental domain $D$ of $\Gamma_P$-action on
$P_0$ has finite diameter $R$.

\hfill

{\bf  Step 3:} 
Fix a non-cusp point $x\in \Abs P$, and fix $z\in D$
(for another choice of $z\in P$ replace $R$ by $R+d(z,D)$).
Let $l_+\subset P$ be a geodesic ray with an end in $x$.
We claim that $l_+\cap P_0$ contains points arbitrarily close to $x$.
The horoballs removed from $P$ to obtain $P_0$ are pairwise disjoint (\ref{_thick_thin_Theorem_}).
If a terminal subray of $l_+$ avoided $P_0$, then, being connected, it would lie in a single horoball $B$.
Its endpoint $x$ would then be the point at which $B$ is tangent to $\Abs$, i.e. a cusp point,
contrary to our assumption. Choose $y_i\in l_+\cap P_0$ with $y_i\to x$.
Since $D$ is a fundamental domain of $\Gamma_P$ acting on $P_0$, there are $\gamma_i\in\Gamma_P$
with $d(y_i,\gamma_i z)\le R$. By \ref{_geode_distance_Claim_} (ii),
$\gamma_i z\to x$, so $x$ belongs to the limit set of $\Gamma_P$.
\endproof

\hfill

\remark
Note that the cusp points of ${\cal H}$ sometimes 
(but not always) also belong to the
limit set of $\Gamma_P$.

\subsection{Cusp points outside of the limit set}

Let  ${\cal P}\subset {\cal H}$ be a CHOPB,
and $P\subset {\Bbb H}^n$ a connected component of its preimage
in the hyperbolic space, equipped with $\Gamma_P=\pi_1({\cal P})$-action.
Denote by $\Lambda\subset \Abs P$ the limit set of $\Gamma_P$-action
on $P$. In Subsection \ref{_abs_of_polyhe_Subsection_},
we have shown that $S:= \Abs P\backslash \Lambda$ is a collection
of cusp points. We don't have much control over $S$;
the following example implies that $S$ can be infinite.

\hfill

\proposition\label{_S_infinite_Theorem_}
Let ${\cal H}$ be a punctured 2-torus
equipped with a complete hyperbolic metric,
and $R\subset {\cal H}$ a simple geodesic connecting
a cusp point to itself. Then ${\cal P}:={\cal H}\backslash R$
is a CHOPB with $\pi_1({\cal P}) =\Z$.
Moreover, if $P\subset{\Bbb H}^2$ is a connected component of the preimage
of ${\cal P}$, then $\Abs P$ contains infinitely many cusp points, while the
limit set of $\pi_1({\cal P})$ acting on $P$ consists of only two points.

\hfill 

\proof
Let $P\subset {\Bbb H}^2$ be a connected component
of the preimage of ${\cal P}$ in the universal cover ${\Bbb H}^2$ of ${\cal H}$.
To avoid confusing the two symbols, note that $P$ is a convex, hence contractible,
subset of ${\Bbb H}^2$, and ${\cal P}= P/\pi_1({\cal P})$ is its quotient
by the stabilizer of $P$; in particular, $P$ is the universal cover of ${\cal P}$.
Cutting the punctured torus ${\cal H}$ along $R$, we see that
${\cal P}$ is diffeomorphic to a cylinder, hence
$\pi_1({\cal P}) =\Z$. The corresponding closed geodesic $c$ in ${\cal P}$
has finite length, hence it corresponds to a hyperbolic isometry
of ${\Bbb H}^2$ preserving $P$. This implies that $ \pi_1({\cal P}) $
has only two limit points, which are the ends of the lift of $c$.
The group $\pi_1({\cal P})$ acts freely on the set of
cusp points of $P$, which all lie in $S$, hence $S$ is infinite.
\endproof


\section{The Apollonian carpet}
\label{_carpet_Section_}


\subsection{Real analytic intervals in the limit set}

\definition
Let $V_\Z$ be a lattice
equipped with a scalar product of signature 
$(1,n)$, and $\Gamma \subset SO^+(V)$
a subgroup commensurable with $SO(V_\Z)$.
Consider a convex hyperbolic orbifold 
${\cal P}\subset {\Bbb H}/\Gamma$ with
polyhedral boundary. Let $P$ be a connected
component of the preimage of ${\cal P}$
in ${\Bbb P} V^+={\Bbb H}$, and
$\Abs P = \bar{P}\cap \Abs$. {\bf The Apollonian carpet} associated to ${\cal P}$ 
is the union of all positive-dimensional
real analytic subvarieties in $\Abs P$.

\hfill

We start by proving that the
Apollonian carpet is a union of geodesic spheres
in $\Abs$; this result follows from a theorem of N. Shah.

\hfill

\theorem\label{_carpet_union_Theorem_}
Let ${\cal H}$ be an arithmetic hyperbolic
orbifold, ${\cal H}= {\Bbb P}V^+/\Gamma$, where
$\Gamma$ is commensurable with 
$SO(V_\Z)$,
and  ${\cal P}\subset {\cal H}$ a convex hyperbolic
orbifold with polyhedral boundary. Denote by $P\subset
\Bbb H$ a connected component of its preimage, and
let $\Gamma_P\subset \Gamma$ be its stabilizer.
Then the Apollonian carpet is a union of spheres
$\Abs \cap {\Bbb P} W_i$, where 
$W_i \subset V$ is a rational 
subspace of signature $(1, k)$.
Moreover, the group $\Gamma_P$
acts on the set $\{W_i\}$ with finitely many orbits.

\hfill

\proof See Subsection \ref{_carpet_proof_Subsection_}. 
\endproof

\hfill

\definition
A {\bf component} of the Apollonian carpet $A$
is a sphere $\Abs \cap {\Bbb P} W_1 \subset A$
which is not contained in another sphere
$\Abs \cap {\Bbb P} W_2 \subset A$ of greater dimension.

\hfill

\remark
From \ref{_carpet_union_Theorem_} it follows that
the Apollonian carpet is the union of its components.

\hfill

\remark  Note that the components of the Apollonian carpet
can be disjoint or tangent, or can intersect transversally;
see examples in Subsection \ref{sec_packing}.

\subsection{Shah's theorem on ergodic properties of the geodesic flow}

Let $ST {\Bbb H}^n$ be the manifold
of unit tangent vectors to the hyperbolic space. We use 
$\Vis:\; ST {\Bbb H}^n\arrow \Abs$ to denote 
the map taking a tangent vector to the limit
point of the corresponding geodesic. If ${\cal H}={\Bbb H}^n/\Gamma$ is a hyperbolic manifold and $(x,v)\in ST {\cal H}$, we 
define $\Vis(x,v)$ by fixing a point in the preimage of $x$ in ${\Bbb H}^n$ and lifting the geodesic to ${\Bbb H}^n$, so that $\Vis$ is defined up to $\Gamma$-action. The definition also makes sense for orbifolds for a suitable definition of $ST$ at the singular points. In our situation every orbifold has a finite cover which is a manifold (\ref{selberg}), and for our purposes we may assume it is already a manifold.

\hfill

\definition
Let $(x, v)\in ST {\Bbb H}^n$,
with $v\in T_x {\Bbb H}^n$ being a unit tangent vector.
Denote by $l_{x,v}\colon [0,+\infty)\to {\Bbb H}^n$ the geodesic starting
at $x$ and tangent to $v$. The
{\bf geodesic flow} is a measure-preserving
flow of diffeomorphisms $\tilde g_t:\; ST {\Bbb H}^n \arrow ST {\Bbb H}^n$
taking $(x, v)$ to the tangent vector $l_{x,v}'(t)$. Similarly, one defines the geodesic flow $g_t$ on ${\cal H}$.

\hfill

\remark By a theorem of E. Hopf, the geodesic flow on a 
hyperbolic manifold is ergodic. The closure of
a given geodesic on a hyperbolic manifold or orbifold might be very pathological, such as
a product of a Cantor set and an interval. 

By Ratner's orbit closure theorem (\cite{_Ratner:Raghunathan_}, see also \cite{_Shah:closures_}), this pathology does not happen when instead of a geodesic, one considers the image in 
${\cal H}$ of a hyperbolic subspace of dimension $k>1$ in ${\Bbb H}^n$: the closure is a closed totally geodesic submanifold of finite volume, which in the arithmetic case is the image of ${\Bbb P}U^+$ for a rational subspace $U$ (\ref{_Ratner_closure_Remark_}). The following theorem by N. Shah goes further in this direction.

\hfill

\theorem  \cite[Theorem 1.2]{_Shah:limiting_}
Let ${\cal H}={\Bbb H}^n/\Gamma$ be a hyperbolic manifold, and 
$C\subset ST {\cal H}$ a real analytic curve, such that $\Vis(C)$
is not a singleton. Then there is a hyperbolic subspace ${\Bbb H}^m\subset {\Bbb H}^n$ with $\Vis(C)\subset \Abs {\Bbb H}^m$
such that the image ${\cal H}_1$ of ${\Bbb H}^m$ in ${\cal H}$ is a closed submanifold, and for any function $f\in C_c(ST{\cal H})$
\begin{equation}\label{eqn_Shah}
	\lim_{t\to \infty}\int_{[0,1]} f(g_t\psi(s))ds=\int_{ST{\cal H}_1}f|_{ST{\cal H}_1} dv,
\end{equation}
where $\psi: [0,1]\to C$ is a parametrization of $C$ by the unit interval and $dv$ the normalized Riemannian volume form. \endproof

\hfill

In particular, taking $f$ which is zero outside of a small neighbourhood of some point $x\in ST{\cal H}_1$, one obtains that $x$ is a limit point of geodesic rays starting at points of $C$, cf. \cite[Remark 1.1]{_Shah:limiting_}. In the same remark, it is explained how to obtain ${\cal H}_1$: choose the smallest subsphere of $\Abs {\Bbb H}^n $ containing $\Vis(C)$, it is the absolute of a hyperbolic subspace ${\Bbb H}^m$, and ${\cal H}_1$ is the closure of the image of  ${\Bbb H}^m$ in ${\cal H}$.



\subsection{Apollonian carpet is a union of geodesic spheres}
\label{_carpet_proof_Subsection_}

We are now ready to prove \ref{_carpet_union_Theorem_}.

\hfill

\proof
As was explained in the discussion above and in \ref{selberg},
we may assume that ${\cal H}$ is a hyperbolic manifold.
Let ${\cal S}_1, \ldots, {\cal S}_m\subset {\cal H}$ be the walls of $\cal P$,
see \ref{_polyhe_orbifo_Definition_}.

\hfill

{\bf Step 1:}
Let $Z\subset\Abs P$ be a real analytic subvariety of positive dimension.
Every point of $Z$ lies on a non-constant real analytic arc contained in $Z$.
Let $c:[0,1]\to\Abs P$ be such an arc, and let $C$ be the image of $c$.
We claim that there exists a hyperbolic subspace $H_1\subset\bar P$
whose image in $\mathcal{H}$ is a closed totally geodesic submanifold
of finite volume and such that $C$ is contained in the absolute of $H_1$.

Fix $z\in P$ and let $\tilde\psi(s)\in ST_z\mathbb{H}$ be the unit vector with
$\mathrm{Vis}(\tilde\psi(s))=c(s)$; let $C_1 = \Im(\tilde\psi)$.
The map $\tilde\psi$ is real analytic, and $\mathrm{Vis}(C_1)=C$ is not
a singleton. Therefore we may let $\psi$ be the composition of
$\tilde{\psi}$ and the projection to $ST{\cal H}$ and apply
\cite[Theorem 1.2]{_Shah:limiting_} recalled above.
We get a hyperbolic subspace $H_1\subset\mathbb{H}$
such that $C$ is contained in the absolute of $H_1$.
It remains to show that we have $H_1\subset\bar P$.
Let ${\cal H}_1$ be the image of $H_1$ in $\cal H$.

Since $P$ is convex and $c(s)\in\bar P$, the geodesic ray
$\{\tilde g_t\tilde\psi(s)\ |\ t\geq 0\}$ lies in $P$; in particular it
never crosses a wall of $P$. Let us show that
${\cal H}_1$ cannot have transversal intersections with any of
the walls ${\cal S}_j$. To this end, define $U\subset ST{\cal H}$ to be
the set of vectors $v$ such that the geodesic segment
$\{g_tv\ |\ 0<t<1\}$ crosses one of the walls ${\cal S}_j$ transversally.
Note that $U$ is open and $g_t\psi(s)\notin U$ for all $t\geq 0$, $s\in[0,1]$.
Applying equation (\ref{eqn_Shah}) from Shah's theorem to
test functions $f\geq 0$ supported in $U$, we conclude that
$ST{\cal H}_1 \cap U = \emptyset$. Therefore ${\cal H}_1$ cannot
have transversal intersections with any of the walls ${\cal S}_j$,
and for any ${\cal S}_j$ the submanifold ${\cal H}_1$ is
either contained in ${\cal S}_j$ or does not intersect it. Finally, note
that $C$ lies both in $\Abs P$ and $\Abs H_1$; hence
$H_1\subset\bar P$. Thus $S:=\Abs H_1$ is a sphere with
$C\subset S\subset\Abs P$.

\hfill

{\bf Step 2:}
Consider the vector subspace $W\subset V$ such that $S= {\Bbb P} W \cap \Abs$. The convex hull of $S$ is 
the hyperbolic subspace ${\Bbb P} W^+\subset {\Bbb P} V^+$.
Note that $\dim W\geq 3$, because $S$ contains the curve $C$.
By \ref{_Ratner_closure_Remark_}, which relies on Ratner's theorem
\cite{_Ratner:Raghunathan_} in the form of \cite{_Shah:closures_},
its image in ${\cal H}= {\Bbb P}V^+/\Gamma$
is closed if and only if 
$\Gamma_W:= \{\gamma\in\Gamma\ \ |\ \ \gamma(W)=W\}$ is a
lattice in $O(W)$, which happens if and only if $W$ is rational. Since this image is ${\cal H}_1$,
which is closed, $W$ must be rational. 

\hfill

{\bf  Step 3:} It remains to show that $\Gamma_P$ acts on the
set of components of the Apollonian carpet with finitely many orbits.
The argument here is analogous to \cite[Section 4]{_AV:Kaw_Mor_};
we only recall it for the sake of completeness.
By Steps 1 and 2, the components are of the form $S= {\Bbb P} W \cap \Abs$,
where $W\subset V$ is a rational subspace, maximal with respect to inclusion.
Denote by $H_S\subset {\cal H}$ the image of ${\Bbb P} W^+$ in ${\cal H}$.
We claim that there is only a finite number of such
submanifolds $H_S\subset {\cal H}$.
By theorems of Mozes--Shah and Dani--Margulis
\cite[Corollaries 1.1, 1.3, 1.4]{_Mozes_Shah_},
the closure of the union of countably many such subsets
is a finite union of subsets of the same shape
(i.e. images of ${\Bbb P}U^+$ for possibly larger rational
$U$). If there is an infinite number of $H_S$, we thus
obtain a closed hyperbolic submanifold ${\cal H}'$
as the closure of some union of them. Such ${\cal H}'$
is different from all the $H_S$.
But the $H_S$ are contained in $\bar{\cal P}$, hence 
${\cal H}'$ is also contained there. It follows that there
exists a subspace $W'\subset V$ such that ${\Bbb P}W'^+$
is contained in the closure of $P$ and ${\cal H}'$ is
the image of ${\Bbb P}W'^+$. Then $\bbP(W')\cap \Abs$ is a new component
of the Apollonian carpet, a contradiction.
\endproof

\subsection{Closure of the Apollonian carpet}

Recall that a group $\Gamma\subset SO^+(1,n)$
is called {\bf elementary} if its limit set has at most 2
points.

\hfill

\theorem\label{thm_cano}
Let $\Gamma\subset SO^+(1,n)$ be a non-elementary
subgroup, and $\Lambda\subset \Abs$ its limit set.
Then the closure of an orbit $\Gamma\cdot x$ 
for any $x\in \Lambda$ is the whole of $\Lambda$.

\proof \cite[Theorem 1.2.24]{_Cano_Mavarrete_Seade_}.
\endproof

\hfill

This result immediately implies the following corollary.

\hfill

\corollary\label{_Apollonian_closure_everything_Corollary_}
Let ${\cal H}$ be a hyperbolic orbifold,
and ${\cal P} \subset {\cal H}$ a CHOPB.
Denote by $\Lambda({\cal P})\subset \Abs$
the limit set of $\Gamma_P=\pi_1({\cal P})\subset \Iso({\Bbb H}^n)$,
and $A\subset \Lambda$ its Apollonian carpet.
Assume that $A$ is non-empty.
Then $\Lambda$ is the closure of $A$.
\endproof


\section{Apollonian carpet of a hyperk\"ahler manifold}
\label{_Ap_carpet_Section_}


\subsection{Apollonian carpet and MBM classes}

\definition
{\bf The Apollonian carpet of a hyperk\"ahler  manifold}
is the union of all geodesic spheres of positive dimension which
belong to the isotropic boundary of its ample cone.

\hfill

\remark
Clearly, the Apollonian carpet of a hyperk\"ahler  manifold
is the Apollonian carpet of the corresponding
convex hyperbolic orbifold with
polyhedral boundary ${\cal P}\subset {\cal H}$,
obtained as the image of ${\Bbb P}\AA_M$ in
${\cal H}= {\Bbb P}\Pos_M/\Mon_I$
(\ref{_kahler_convex_hy_polyhe_Theorem_}).

\hfill

The Apollonian carpet of a hyperk\"ahler  manifold
can be described explicitly in terms of its BBF form
and MBM classes, as follows.

\hfill

\theorem\label{_Apollonian_carpet_explicitly_Theorem_}
Let $M$ be a projective hyperk\"ahler manifold, and
$\NS(M)= H^{1,1}(M, \R) \cap H^2(M, \Z)$ its 
N\'eron--Severi lattice. Consider the set
of all maximal rational subspaces $W_i \subset \NS_\R$
of signature $(1,k)$, $k\geq 2$ such that for any
MBM class $\eta \in \NS(M)$ with $\eta\not\!\!\bot W_i$, 
we have that $W_i^\bot +\eta$ is not negative definite. Then 
the sphere $S_i= \Abs \cap {\Bbb P}W_i$ 
 lies in the closure of some ample chamber. By \ref{chambers}, this chamber is the ample cone of a birational model of $M$ or the image of such an ample cone under an element of $\Mon_I$; in both cases it can be seen as the ample cone of $(M, I_i)$ for a complex structure $I_i$.
Moreover $S_i$ belongs to the Apollonian carpet of 
$(M,I_i)$, and  all components of the
Apollonian carpet of $(M,I_i)$ are obtained this way. Finally, the group of automorphisms of $(M,I_i)$ contains a lattice in 
$O(1, k)=\Aut W_i$, which is commensurable with $O(W_i, \Z)$.

\hfill

\remark\label{face} The subspace $W_i$ as above can lie in $\eta^{\bot}$ for some MBM class $\eta$. The space generated by all such $\eta$ for a given  $W_i$ is negative definite. If the addition of an extra MBM class $\eta'$ destroys negativity, this means that $W_i$ supports a face of an ample chamber which does not contain other faces, in other words, $W_i^+$ itself is a minimal face of the ample chamber. Changing the complex structure if necessary, we may assume that the chamber is the ample cone.

\hfill

{\bf Proof of the theorem, step 1:}

When $W_i^+$ is a face of the ample cone as in \ref{face}, clearly $\Abs {\Bbb P}W_i^+$ is a sphere on the boundary of the ample cone and a component of the Apollonian carpet. 

In general, let $U_1, U_2 \subset V$ be negative definite subspaces
of a space of signature $(1,n)$. Then ${\Bbb P}(U_1^\bot)^+$
intersects ${\Bbb P}(U_2^\bot)^+$ if and only if
$U_1+U_2$ is negative definite. If $U_1+U_2$
is degenerate, the spheres $\Abs {\Bbb P}(U_1^\bot)^+$
and $\Abs {\Bbb P}(U_2^\bot)^+$ are tangent.

 Take a $W_i$ which is not contained in  $\eta^\bot$ with $\eta$ MBM. The condition of the theorem: $W_i^\bot+ \eta$ is not negative definite means that $S_i$ does 
not meet any orthogonal hyperplane to an MBM class except maybe tangentially. This means that $S_i$ is contained in a chamber which, as before, we assume to be the ample cone. Therefore $S_i$ is again a component of the Apollonian carpet.

{\bf Step 2:}

Let $S_i$ be a component of the Apollonian
carpet of $(M,I)$. Then $S_i=\Abs \cap {\Bbb P}W_i$,
where $W_i$ is a rational subspace. Since  $W_i$ intersects the positive cone, it must have signature $(1,k)$.
Also, $S_i$ is contained in the closure of the ample cone. The ample cone is cut inside the positive cone by the orthogonal 
hyperplanes to MBM classes. We derive that for $\eta$ MBM, $S_i$ cannot intersect $\eta^{\bot}$ transversally (it is disjoint, tangent to or contained in $\eta^{\bot}$). This amounts to saying  that  $W_i^\bot +\eta$ is not negative definite when $W_i$ is not orthogonal to $\eta$. 

{\bf Step 3:}

It remains to prove the last assertion of the theorem. Replacing $I$ by $I_i$, we may assume that $S_i$ is a component of the Apollonian carpet of $(M,I)$, that is, $S_i$ lies in the closure of ${\Bbb P}\AA_M$. We first recall the setting of \ref{_kahler_convex_hy_polyhe_Theorem_}.

Denote by $\Gamma_I$ the image of the Hodge monodromy group $\Mon_I$ in
the orthogonal group $O(H^{1,1}(M,\Z))$. By \ref{_Mon_I_arithm_Proposition_}, the map $\Mon_I\to\Gamma_I$ has finite kernel and $\Gamma_I$ has finite index in $O(H^{1,1}(M,\Z))$, so $\Gamma_I$ is an arithmetic lattice acting on ${\Bbb P}\Pos_M$, and
${\cal H}:= {\Bbb P} \Pos_M/\Gamma_I$
is an arithmetic hyperbolic orbifold (the monodromy orbifold of \ref{def_ample_poly}).
The hyperplanes orthogonal to the MBM classes of type $(1,1)$ decompose ${\Bbb P}\Pos_M$ into ample chambers, one of which is ${\Bbb P}\AA_M$, and $\Gamma_I$ permutes these chambers (\ref{chambers}). Let $\Gamma\subset\Gamma_I$ be the stabilizer of the ample cone $\AA_M$. By \ref{_kahler_convex_hy_polyhe_Theorem_}, the image ${\cal P}$ of ${\Bbb P}\AA_M$ in ${\cal H}$ is a convex hyperbolic orbifold with polyhedral boundary, ${\Bbb P}\AA_M$ is its universal cover, and $\pi_1({\cal P})=\Gamma$. By \ref{chambers}, the elements of $\Gamma_I$ which lift to automorphisms of $(M,I)$ are exactly the elements of $\Gamma$; since the kernel of $\Aut(M,I)\to \Gamma$ is finite (\ref{_kahler_convex_hy_polyhe_Theorem_}), it suffices to show that $\Gamma$ contains a subgroup of finite index in the stabilizer of $W_i$ in $\Gamma_I$, and that this stabilizer is a lattice in $O(W_i)$ commensurable with $O(W_i,\Z)$.

Let $S_i$ be a component of the Apollonian carpet and $W_i$ the associated linear subspace as above. Set $\Gamma_{W_i}:= \{ \gamma \in \Gamma_I\ \ |\ \ \gamma(W_i)= W_i\}$.
Since $W_i$ is rational and $\Gamma_I$ is arithmetic, the restriction of $\Gamma_{W_i}$ to $W_i$ is a lattice in $O(W_i)$, commensurable with $O(W_i,\Z)$ (by the Borel--Harish-Chandra theorem, see \cite{_Morris:Ratner_}); the kernel of the restriction map is finite, because it acts faithfully on the negative definite lattice $W_i^\bot\cap H^{1,1}(M,\Z)$. If $S_i$ is not contained in the orthogonal to an MBM class, then $S_i$ intersects all
$\eta^{\bot}$, where $\eta$ is MBM, at most tangentially, and therefore an element $\gamma$ of $\Gamma_{W_i}$, preserving $S_i$, must also 
preserve its chamber, i.e. the ample cone. Thus $\gamma$ lifts to an automorphism of $(M, I)$. If $S_i$ is contained in the orthogonal 
to an MBM class, then $\P \Pos(W_i)$ is a face of the ample cone, and  $\gamma\in \Gamma_{W_i}$ can also take the ample cone to another chamber adjacent to $\P \Pos(W_i)$, but these are finitely many. So a finite index subgroup of $\Gamma_{W_i}$ lifts to a subgroup in 
$\Aut(M,I)$. \endproof

\subsection{The Baragar gasket}

In the proof of \ref{_Apollonian_carpet_explicitly_Theorem_}
 we have seen that the Apollonian carpet has two types of components: the ones which come from the faces of the ample cone, and the ones which do not. 

More precisely, let $F$, $\dim F > 2$, be a
face of the ample cone which is minimal with respect to inclusion,
that is, $F=L^+$ for a vector subspace $L$ that is the intersection
of orthogonal complements to MBM classes and such that $F$ does not meet other orthogonal complements
to MBM classes except possibly tangentially.
Then ${\Bbb P}L^+\cap \Abs$ is a component of the Apollonian carpet.

In the case of K3 surfaces, the union of such components has been studied by A. Baragar, see e.g. \cite{_Baragar1_}. Hence the following definition.

\hfill

\definition\label{def_Baragar}
The union of all components of the Apollonian carpet
obtained this way is called {\bf the Baragar
  gasket}. 

\hfill

\remark\label{rem_ghost}
By definition, the Baragar gasket is contained in the Apollonian carpet,
but in general it is a proper subset. The complement consists of the components
${\Bbb P}W\cap\Abs$ of the carpet (\ref{_Apollonian_carpet_explicitly_Theorem_})
for which $W$ is not an intersection of orthogonal complements to MBM classes, that is,
of spheres which lie on the isotropic boundary of the ample cone but are not
of the form ${\Bbb P}L^+\cap\Abs$ for a minimal face $L^+$ as above. The complement may well be non-empty: see \ref{nonbaragar}; in
\ref{intersect} two non-Baragar components even intersect transversally;
and in Example 1 of Subsection \ref{sec_illustr} the Apollonian carpet is
non-empty while the Baragar gasket is empty.
This phenomenon has already been observed in the context of circle packings,
where the components of the complement have appeared in the literature under
the name of {\bf ghost circles} (\cite{St}).

\hfill

\remark
The Baragar components do not intersect transversally, but
they might be tangent.  Other components of the
Apollonian carpet can intersect (\ref{intersect}), but they cannot
intersect the Baragar components transversally. In the
CHOPB setting, the Baragar components correspond
to the hyperspheres which are obtained as isotropic
boundaries of faces of the CHOPB which don't contain
faces of smaller dimension, and hence are given
by a hyperbolic subspace which lies fully on
a boundary of the CHOPB.

\hfill

\remark
This gives an effective and easy way to 
determine the Apollonian carpet for particular manifolds
when we have good control over the BBF form, the
N\'eron--Severi lattice and the MBM classes. Note that 
by Nikulin's theorem (\cite{_Nikulin_}),  any even quadratic 
lattice of signature $(1, n)$, $n\leq 10$
can be realized as the N\'eron--Severi lattice for an appropriate K3
surface. Note also that for a K3 surface, MBM classes are integer
classes with square $-2$.

\subsection{Examples}\label{sec_packing}

We give several examples illustrating the shape of the Apollonian carpet and Baragar gasket. For simplicity, we restrict ourselves to K3 surfaces.
As we have already observed, a sphere on the boundary of the projectivized positive cone is the intersection of this boundary with a linear subspace $W$, and the BBF form $q$ (the intersection
form in this case) is negative definite on $W^\bot$. This sphere belongs to an ample chamber if it is either 
disjoint from the walls, or intersects them tangentially. This condition means that adding a $(-2)$-class $z$ to $W^\bot$ destroys negativity:
$z^{\bot}$ is disjoint from the sphere if and only if $q$ is of signature $(1, \dim W^\bot)$ on $W^\bot+z$, and tangent to the sphere if and only if $q$ is 
negative semi-definite on $W^\bot+z$. Finally, such a sphere is a Baragar component when $W^\bot$ is generated by integral classes of square $-2$, and not a Baragar component otherwise.

To calculate the signature of $q$ on the relevant subspaces of the real N\'eron--Severi group, we prove some auxiliary results.
Recall that the discriminant group of a non-degenerate quadratic lattice $\Lambda$ is
$D(\Lambda) = \Lambda^\vee / \Lambda$, where $\Lambda^\vee$ is the dual lattice and
the embedding $\Lambda \hookrightarrow \Lambda^\vee$ is induced by the bilinear form;
its order is $|D(\Lambda)| = |\det q|$.

\hfill

\lemma\label{lemma_discr}
Assume that $L\subset \Lambda$ is a primitive non-degenerate sublattice, and denote by $A$ the
quotient $\Lambda/L$. Then we have the following exact triples:
\begin{eqnarray}
&0\to A\to \Lambda^\vee/L \to D(\Lambda)\to 0&\label{eqn_triple1}\\
&0\to A^\vee\to \Lambda^\vee/L \to D(L)\to 0&\label{eqn_triple2}
\end{eqnarray}

\hfill

\begin{proof}
By the definition of $D(\Lambda)$ we have an exact triple
$$
0\to \Lambda\to \Lambda^\vee\to D(\Lambda)\to 0.
$$
Taking the quotient of the first two terms by $L$, we obtain (\ref{eqn_triple1}).

Since the sublattice $L$ is primitive, the quotient $A = \Lambda/L$ is free, so we have an exact triple
$$
0\to A^\vee\to \Lambda^\vee\to L^\vee\to 0.
$$
Taking the quotient of the last two terms by $L$, we obtain (\ref{eqn_triple2}).
\end{proof}

\hfill

\proposition\label{prop_discr}
Let $(\Lambda, q)$ be a quadratic lattice of 
rank $r$ for some $r\ge 3$, and let $p$ be a prime number.
Assume that the $\mathbb{F}_p$-vector
space $D(\Lambda) \otimes \mathbb{F}_p$ has dimension $r-k$. Then for any
non-degenerate sublattice $L\subset\Lambda$ of rank $l>k$, $p$ divides $|D(L)|$.


\hfill

\begin{proof}
First assume that $L$ is primitive. Applying \ref{lemma_discr} and using our assumption
on $D(\Lambda)$, we deduce from the sequence (\ref{eqn_triple1}) by tensoring it with $\mathbb{F}_p$
that the dimension of $(\Lambda^\vee/L)\otimes\mathbb{F}_p$ is at least $r-k$. On the other
hand, $A^\vee$ is a free abelian group of rank $r-l$, and we deduce from the sequence (\ref{eqn_triple2})
by tensoring it with $\mathbb{F}_p$ that $D(L)\otimes\mathbb{F}_p$ is at least one-dimensional.
Hence $D(L)$ contains non-trivial $p$-torsion, and $p$ divides $|D(L)|$.

In general, let $\tilde L = (L\otimes\bbQ)\cap\Lambda$ be the saturation of $L$. It is a primitive
non-degenerate sublattice of rank $l$, so $p$ divides $|D(\tilde L)|$ by the above. Since
$|D(L)| = |\det q|_L| = [\tilde L:L]^2\,|D(\tilde L)|$, $p$ divides $|D(L)|$ as well.
\end{proof}

\hfill

We will use the following elementary observation. Let $V$ be a real vector space with a form
of signature $(1,n)$. Then a degenerate subspace $U\subset V$ is negative semi-definite.
Indeed, its kernel contains a non-zero isotropic vector $u$, and $U\subset u^\perp$, on which
$q$ is negative semi-definite.

\hfill

\corollary\label{rk2} Let $(\Lambda, q)$ be a quadratic lattice of signature $(1,r-1)$ for some $r\ge 3$, with 
$D(\Lambda) \otimes \mathbb{F}_p$ of dimension $r-1$. Let $L$ be a rank-two 
sublattice generated by $q$-negative vectors $v_1$, $v_2$ with $q(v_1)q(v_2)<p$. Then $L$ is either degenerate and negative semi-definite, 
or has signature $(1,1)$.

\hfill

\begin{proof}
If $L$ is degenerate, it is negative semi-definite by the observation above. Otherwise, since $L$ contains the negative
vector $v_1$, it is either of signature $(1,1)$ or negative definite. In the latter case, by the Cauchy--Schwarz inequality,
$$
0 < |D(L)| = q(v_1)q(v_2) - q(v_1,v_2)^2 \le q(v_1)q(v_2) < p .
$$
On the other hand, \ref{prop_discr} with $k=1<2$ gives $p\mid |D(L)|$, a contradiction.
\end{proof}

\hfill

\corollary\label{rk3} Let $(\Lambda, q)$ be a quadratic lattice of signature $(1,r-1)$ for some $r\ge 3$, with
$D(\Lambda) \otimes \mathbb{F}_p$ of dimension $r-2$, where $p\geq 5$. Let $L$ be a rank-three sublattice generated by
three elements $v_1,v_2,v_3$ of square $-2$ with $q(v_1,v_2)=0$. Then $L$ is of signature $(1,2)$ or degenerate negative
semi-definite of signature $(0,2)$.

\hfill

\begin{proof}
The plane $\langle v_1,v_2\rangle$ is negative definite. So if $L$ is degenerate, it is negative semi-definite with
a one-dimensional kernel, by the observation above. Otherwise $L$ has signature $(1,2)$ or is negative definite. Assume the latter,
and put $c=q(v_1,v_3)$, $d=q(v_2,v_3)$. Then
$$
|D(L)| = \det\begin{pmatrix} 2&0&-c\\ 0&2&-d\\ -c&-d&2\end{pmatrix} = 8-2(c^2+d^2).
$$
Negative definiteness gives $c^2+d^2\le 3$, hence $c^2+d^2\in\{0,1,2\}$ and $|D(L)|\in\{4,6,8\}$.
None of these is divisible by a prime $p\ge 5$. This contradicts \ref{prop_discr} with $k=2<3$.
\end{proof}

\hfill

\remark The assumptions $p\ge 5$ and $v_1\perp v_2$ are exactly what is needed to exclude
$p\mid|D(L)|$ in the computation above.

\hfill

\example\label{example_first} Consider the lattice $\Lambda$ of rank four with the quadratic form $q(x,y,z,w) = 2px^2 - 2y^2 - 2pz^2 - 2pw^2$
for some prime number $p\geq 5$. This is an even lattice of signature $(1,3)$. By the classical results
about lattice embeddings (\cite{_Nikulin_}, see also \cite[Section 3.3]{_AV:automorphisms_}), $\Lambda$ can be primitively embedded into the K3 lattice $\Lambda_{K3}$,
therefore $\Lambda$ can be realized as the N\'eron--Severi lattice of an algebraic K3 surface. Set $V=\Lambda\otimes \R$ and consider ${\Bbb H}^3={\Bbb P} V^+$. Then $\Abs=S^2$. 
 It is clear that the set of classes of square $-2$ is non-empty. By \ref{rk2}, the circles on $\Abs$ cut out by the orthogonal hyperplanes to $(-2)$-classes 
 $v_1$ and  $v_2$ do not intersect, except possibly tangentially. 
Therefore we obtain a circle packing on the boundary of the positive cone of the K3 surface
and a non-trivial Baragar gasket (see \ref{def_Baragar}) in the closure of its ample cone.

\hfill

\example\label{example_lat} Note that in the example above, some spheres of the Baragar gasket may be tangent.
Let us construct an example where all of the spheres are disjoint. Note that the spheres
defined by two vectors $v_1, v_2\in \Lambda$ (as $v_i^\bot\cap \Abs$) are tangent if and only if the restriction
of $q$ to $\langle v_1, v_2\rangle$ is degenerate. In particular, this implies that $\Lambda$
contains non-trivial isotropic vectors (they correspond to the points of tangency of the spheres).
Therefore to construct our example it is sufficient to produce an even lattice without non-zero
isotropic elements such that it contains classes of square $-2$ and satisfies the conditions
of \ref{prop_discr} for some prime $p\ge 5$.

Start from the rank-four lattice $\Lambda_0$ with the quadratic form $q_0(x,y,z,w) = 3x^2 - y^2 -7z^2-7w^2$.
Let us observe that $\Lambda_0$ does not contain non-zero isotropic elements. To see this, assume
that $x, y, z, w$ are coprime integers satisfying $3x^2 - y^2 -7z^2-7w^2 =0$. This implies
that $3x^2 \equiv y^2\,(\mathrm{mod}\,\, 7)$ and since $3$ is not a square mod $7$, we deduce that
$x \equiv y \equiv 0\,(\mathrm{mod}\,\, 7)$. Writing $x = 7x_0$ and $y = 7y_0$ we see that
$21x_0^2 - 7y_0^2 -z^2-w^2 =0$. It follows that $z^2+w^2\equiv 0\,(\mathrm{mod}\,\, 7)$,
and since $-1$ is not a square mod $7$, we deduce that $z \equiv w \equiv 0\,(\mathrm{mod}\,\, 7)$,
arriving at a contradiction.

Now let $\Lambda$ be the rank-four lattice with the quadratic form $q(x,y,z,w) = 6p^2x^2-2y^2-14p^2z^2-14p^2w^2$,
i.e. $q$ is obtained from $2q_0$ by multiplying the first, third and fourth summands by $p^2$, for
some prime $p\ge 5$. Then $(\Lambda, q)$ is clearly a sublattice of $(\Lambda_0, 2q_0)$, therefore
it does not contain non-trivial isotropic elements. It is also an even lattice of signature $(1,3)$
representing $-2$ and satisfying the assumptions of \ref{prop_discr} for $k=1$. By the same argument
as above, this lattice can be realized as the N\'eron--Severi lattice of an algebraic K3 surface,
giving an example of the Baragar gasket with all spheres disjoint. 

\hfill

\example\label{nonbaragar} Let us now give an example of an Apollonian carpet which also has non-Baragar components. Let 
$\Lambda$ be the rank-four lattice with the quadratic form $q(x,y,z,w) = 46x^2 - 2y^2 - 46z^2 - 46w^2$ (i.e. the form from \ref{example_first}
with $p=23$). This form obviously represents both $-2$ and $-4$. Moreover, by \ref{rk2}, if $q(v_1)=-2$ and $q(v_2)=-4$, the circles cut out on the absolute by $v_1^{\bot}$ and $v_2^{\bot}$ do not intersect except maybe tangentially. This means that the latter circle lies entirely on the
isotropic boundary of the projectivization of the ample cone, i.e. is a component of the Apollonian carpet. Notice that this is not a Baragar component, because $q(v_2)=-4$ and all MBM classes are of square $-2$.

\hfill

\example\label{codim} The components of the Apollonian carpet can be of codimension more than one. Consider a K3 surface with Picard lattice $\Lambda$ of rank 5, with 
$q(x,y,z,w,u)=4px^2-2y^2-2z^2-4pw^2-4pu^2$, where $p\ge 5$. Here $e_2$ and $e_3$ are orthogonal,
of square $-2$, and $D(\Lambda)\otimes\mathbb{F}_p$ has dimension $3=r-2$, so by \ref{rk3} the orthogonal to $\langle e_2,e_3 \rangle$ does not intersect the other walls of the ample chambers except possibly tangentially, therefore gives a codimension-two Baragar component of the Apollonian carpet.

\hfill

\example\label{intersect} The non-Baragar components can also intersect transversally. Take the rank-four lattice
$\Lambda$ with the quadratic form $2q$, where $q=5x^2-y^2-5z^2-5w^2$. The $(-2)$-classes of $(\Lambda,2q)$ are the vectors
$z$ with $q(z)=-1$. Consider $v_1=(1,2,1,0)$ and $v_2=(1,2,0,1)$. Here $q(v_1)=q(v_2)=-4$ and $q(v_1,v_2)=1$, so
the restriction of $q$ to $\langle v_1,v_2\rangle$ is negative definite. Hence the circles
$S_i = v_i^\perp\cap\Abs$ meet transversally, in two points.

The $\mathbb F_5$-space $D(\Lambda,q)\otimes\mathbb F_5$ has dimension $3=r-1$, and $q(v_i)q(z)=4<5$ for every
$(-2)$-class $z$. By \ref{rk2} applied to $q$, the lattice $\langle v_i,z\rangle$ is never negative definite. Therefore
$W_i=v_i^\perp$ satisfies the assumptions of \ref{_Apollonian_carpet_explicitly_Theorem_}, and $S_i$ is a component
of the Apollonian carpet of some ample chamber $C_i$. It is not a Baragar component, because $v_i$ is primitive and
$\Z v_i$ contains no $(-2)$-class.

The two circles lie on the boundary of the same chamber. Indeed, suppose $C_1\neq C_2$. Then they are separated by
some wall $z^\perp$, with $z$ a $(-2)$-class: $q(z,\cdot)>0$ on $C_1$ and $q(z,\cdot)<0$ on $C_2$. Passing to closures,
$q(z,\cdot)\ge 0$ on $S_1$ and $q(z,\cdot)\le 0$ on $S_2$, so $S_1\cap S_2\subset S_z := z^\perp\cap\Abs$. The circle
$S_z$ is different from $S_1$, since $q(z)=-1$ while every multiple of the primitive vector $v_1$ has $q=-4m^2$.
By \ref{rk2}, $S_z$ meets $S_1$ at most tangentially, hence in at most one point. This contradicts
$|S_1\cap S_2|=2$. Thus $C_1=C_2$, and $S_1$, $S_2$ are two transversally intersecting components of the
Apollonian carpet of the same K3 surface. These circles are examples of the ghost
circles of \cite{St} (\ref{rem_ghost}). See Figure \ref{fig_3} in the next subsection.

\hfill

\remark The reason why we have only considered K3 surfaces is that the description of MBM classes is particularly simple in this case:
indeed these are precisely the elements of square $-2$ in the lattice. But it is clear from the discussion above that the same method should yield analogous examples for a hyperk\"ahler manifold as soon as one has a description of its MBM classes. Indeed it is known (\cite{_AV:Kaw_Mor_})  that MBM classes have bounded BBF square, so setting the prime $p$ as above large enough allows one to produce lattices with desired  properties. Strictly speaking, there is an additional difficulty in proving that a given lattice embeds into the second cohomology lattice of a 
hyperk\"ahler manifold, as Nikulin's theorem requires unimodularity of this larger lattice. But there exist some ways to get around this difficulty, see e.g. \cite{_AV:automorphisms_}.

\subsection{Illustrations}\label{sec_illustr}

In this subsection we provide some examples of images of the sphere
packings on the boundary of the hyperbolic 3-space. The source code
of the scripts used to produce the images is available at \cite{Source}.
All of the presented configurations of spheres may appear
in the N\'eron--Severi groups of algebraic K3 surfaces (some of the lattices considered below
are odd; one has to multiply them by two before embedding into the K3 lattice). We hope that
our illustrations will
give an idea of what the isotropic boundary of the ample cone of a hyperk\"ahler manifold
may look like, more vividly than the mere formal description.

We start by describing the procedure that was used to produce all the figures below. We fix
a quadratic lattice $(\Lambda, q)$ of rank $4$ and signature $(1, 3)$. We fix
an integer $d<0$ and let $\cS = \{v\in \Lambda\st q(v) = d \}$. The collection $\cS$
defines an arrangement of hyperplanes $H_v = v^\perp \subset V = \Lambda \otimes \bbR$, where $v\in \cS$.

We fix an isomorphism $\Lambda \simeq \bbZ^4$, identifying $q$ with a symmetric matrix $A$.
We denote by $\nu_0$, $\nu_1$, $\nu_2$ and $\nu_3$ the eigenvectors of $A$ with eigenvalues
$\lambda_0 > 0 > \lambda_1 \ge \lambda_2 \ge \lambda_3$, normalized by the condition $|q(\nu_i)| = 1$.
We identify $V$ with $\bbR^{1,3}$ using the basis $\nu_0,\ldots,\nu_3$,
so that for $x = (x_0, x_1, x_2, x_3)\in \bbR^{1,3}$ we have $q(x) = x_0^2 - x_1^2 - x_2^2 - x_3^2$.

The hyperbolic 3-space $\bbH^3 = \bbP(V^+)$ is identified with the open unit ball in $\bbR^3$
via the stereographic projection from the point $(-1,0,0,0)$. The inverse of this projection is
given by
$$
(y_1, y_2, y_3)\mapsto (1+|y|^2: 2y_1: 2y_2: 2y_3),
$$
where $|y|$ denotes the standard Euclidean norm in $\bbR^3$. Under the stereographic projection the boundary $Q$
of $\bbH^3$ is identified with the unit sphere $S^2\subset \bbR^3$ and we have
$$
Q\cap H_v = \left\{ y\in S^2\st v_0 = \sum_{i=1}^3 v_i y_i \right\},
$$
where $v = (v_0, v_1, v_2, v_3)\in \bbR^{1,3}$.

As we have discussed above, the collection of hyperplanes $H_v$, $v\in \cS$ is locally finite in $\bbH^3$,
and it splits $\bbH^3$ into a countable set of open chambers.
We will study the boundary of the $\cS$-chamber in $\bbH^3$ that contains the origin, which corresponds to the point $(1:0:0:0)\in \bbP(V^+)$.
More precisely, since the origin might lie on one of the walls $H_v$, we will study the boundary
of the union of all $\cS$-chambers that contain the origin in their closure. The number of such chambers
is finite and we denote their union by $\CC_0\subset \bbH^3$. We also define $\cS' = \{v\in\cS\st v_0 > 0 \}$,
so that $\overline{\CC_0} = \cap_{v\in\cS'} H_v^+$. We let $Q_v = H_v^- \cap Q$ for $v\in \cS'$.
Then the isotropic part of the boundary of $\CC_0$
is $$\overline{\CC_0} \cap Q = S^2\setminus \cup_{v\in \cS'} Q_v.$$

We further consider the stereographic projection from the point $(1,0,0)\in S^2$ onto $\bbR^2$.
The inverse of that projection identifies $\bbR^2$ with the complement of $(1,0,0)\in S^2$:
$$
(z_1, z_2) \mapsto (y_1, y_2, y_3) = \left(\frac{|z|^2 - 1}{|z|^2 + 1}, \frac{2z_1}{|z|^2 + 1}, \frac{2z_2}{|z|^2 + 1}\right).
$$
Under this projection $Q_v$ is mapped to the subset
$$
D_v = \{z\in \bbR^2\st |z|^2(v_1 - v_0) + 2v_2z_1 + 2v_3z_2 > v_0 + v_1\},
$$
which is either the interior or the exterior of the circle with centre at the
point $(v_2/(v_0-v_1), v_3/(v_0-v_1))$ of radius $|q(v)|^{1/2} |v_1 - v_0|^{-1}$,
or a half-plane, depending on the non-vanishing and the sign of $v_1-v_0$. Slightly abusing the terminology,
we will call $D_v$ a $v$-disc (although it might not be really a disc, but a ``hole'' in $\bbR^2$ or a half-plane).
The isotropic boundary of $\CC_0$ under our projection is identified with
the complement $\mathcal{B}_0$ to the union of $v$-discs for all $v\in\cS'$.

The $v$-discs are partially ordered by inclusion as subsets of $\bbR^2$, and this defines a partial order
on $\cS'$. It is clear that any increasing chain with respect to this partial order stabilizes (because
of the local finiteness of the hyperplanes $H_v$ in $\bbH^3$). 
Hence to describe the isotropic boundary $\mathcal{B}_0$ it is sufficient to consider only the maximal $v$-discs.
The maximal $v$-discs correspond to the faces of $\overline{\CC_0}$.

Listing all elements of the set $\cS'$, we obtain increasing finite collections of maximal $v$-discs,
and the complements to such finite collections give approximations to $\mathcal{B}_0$. In the pictures
below, we show such finite collections of $v$-discs inside the unit disc in $\bbR^2$ (which
corresponds to the hemisphere in $S^2\subset \bbR^3$ given by the inequality $y_1 < 0$). The complement
to the $v$-discs (black points in the pictures) is the approximation of the isotropic boundary of $\CC_0$.

\subsubsection{Example 0.} \label{_example_Apollonian_Subsubsection_}
The classical Apollonian circle packing is a configuration of circles in the
plane that is obtained iteratively by the following procedure. One starts from a triple of
pairwise tangent circles having three distinct points of tangency.
As the first step one adds the two circles that are tangent to the given three.
Then one repeats the process for all triples of pairwise tangent circles in the configuration obtained
in the previous step. One can obtain the Apollonian configuration starting from the
lattice with the intersection matrix (see \cite{_Baragar1_})

$$
A = \begin{pmatrix}
-1 & 1 & 1 & 1\\
1 & -1 & 1 & 1\\
1 & 1 & -1 & 1\\
1 & 1 & 1 & -1
\end{pmatrix}
$$

\begin{figure}
\centering
\includegraphics[width=8cm]{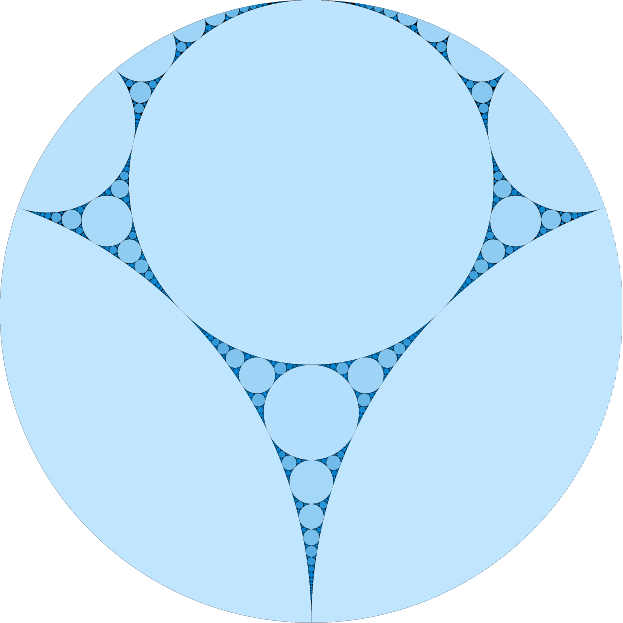}
\caption{The classical Apollonian gasket (Example 0)}\label{fig_0}
\end{figure}

We take $d = -1$. The result is shown in Figure \ref{fig_0}.

\begin{figure}[ht]
\centering
\includegraphics[width=10cm]{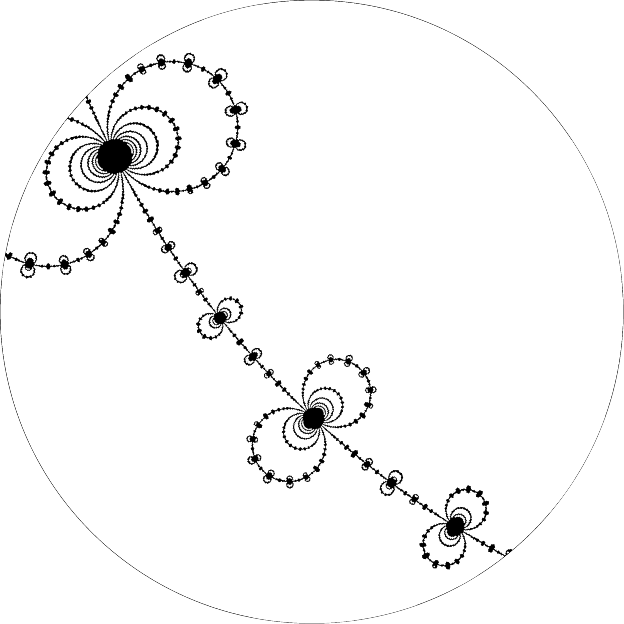}
\caption{The isotropic boundary is a fractal (Example 1)}\label{fig_1}
\end{figure}

\subsubsection{Example 1.} In a typical situation, the maximal $v$-discs may overlap, and the isotropic boundary
(when non-empty) is a fractal. Figure \ref{fig_1} shows such an example, where $d = -1$ and the quadratic form $q$ is
given by the matrix

$$
A = \begin{pmatrix}
-1 & 2 & 0 & 0\\
2 & -1 & 1 & 0\\
0 & 1 & -2 & 1\\
0 & 0 & 1 & -2
\end{pmatrix}
$$

\begin{figure}[h]
\centering
\includegraphics[width=10cm]{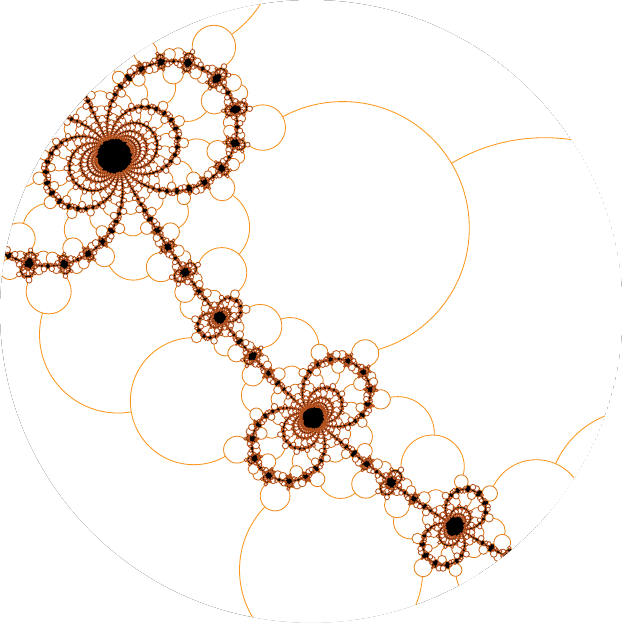}
\caption{The discs that cut out the fractal (Example 1)}\label{fig_2}
\end{figure}

In Figure \ref{fig_2} we show the overlapping $v$-discs that cut out the isotropic boundary.
One may guess from Figure \ref{fig_1} that in this example the isotropic boundary contains infinitely many circles (which form
what we call the Apollonian carpet), but since the $v$-discs overlap, none of their
boundary circles is contained in the isotropic boundary (so the Baragar gasket is empty).

\newpage

\begin{figure}[h]
\centering
\includegraphics[width=10cm]{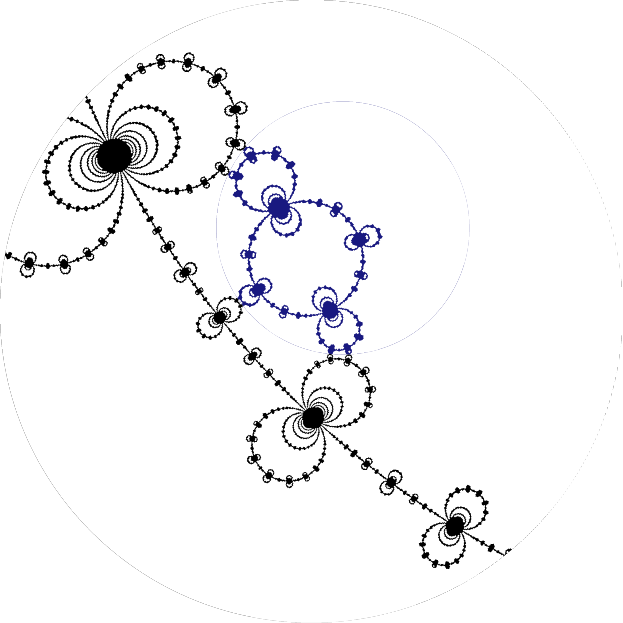}
\caption{The isotropic boundaries of two different chambers (Example 1)}\label{fig_1prime}
\end{figure}

Note that Figure \ref{fig_1} shows (part of) the isotropic boundary of only one chamber.
For the other chambers the isotropic boundary is obtained by a composition of reflections
in some $(-1)$-vectors. Figure \ref{fig_1prime} shows the isotropic boundary for two different chambers.

\subsubsection{Example 2.}

Here is a pictorial representation of \ref{intersect}. 



\newpage

\begin{figure}[h]
\centering
\includegraphics[width=10cm]{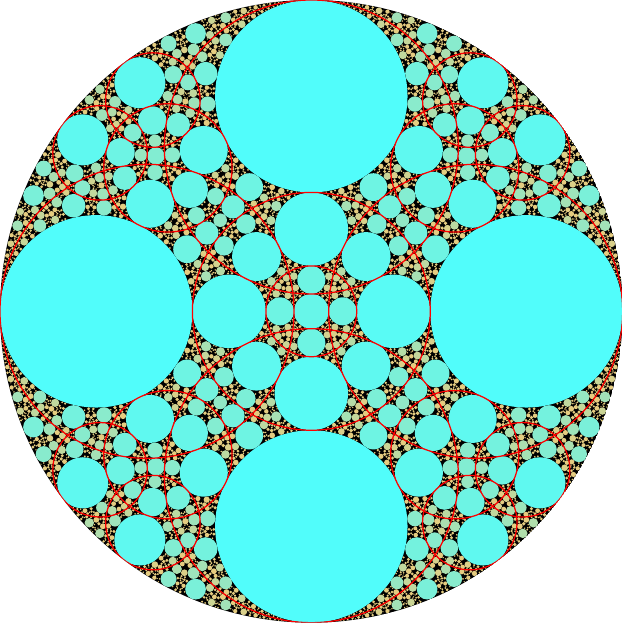}
\caption{A configuration with intersecting ghost circles (Example 2)}\label{fig_3}
\end{figure}

The maximal $v$-discs are disjoint, although their
boundary circles may be tangent. These circles form the Baragar gasket. The circles $C_w = Q\cap H_w$ for primitive elements $w\in \Lambda$ with $q(w) = -4$ are part of the Apollonian carpet, but not of the Baragar gasket, and some of them intersect transversally (this is shown in red in Figure \ref{fig_3}).
Following \cite{St}, we call these circles {\bf ghost circles} (see \ref{rem_ghost}).


\newpage

\begin{figure}[h]
\centering
\includegraphics[width=10cm]{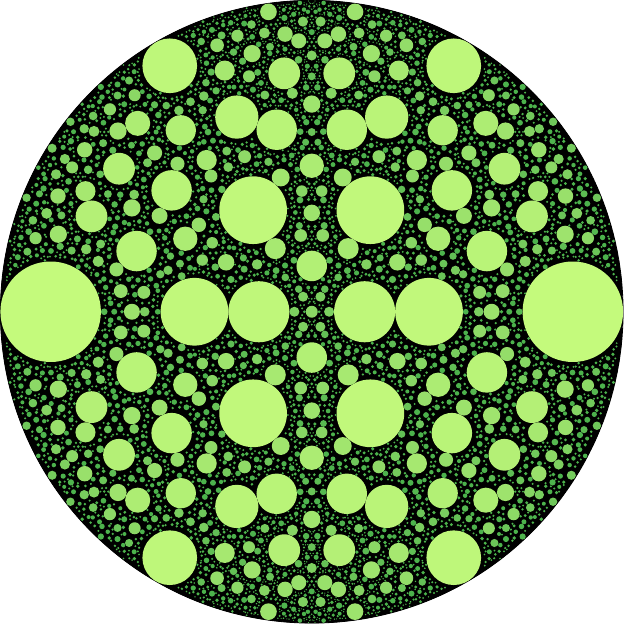}
\caption{An example of a packing (Example 3)}\label{fig_4}
\end{figure}

\subsubsection{Example 3.}

The quadratic form is given by

$$
A = \begin{pmatrix}
p^3 & 0 & 0 & 0\\
0 & -p^2 & 0 & 0\\
0 & 0 & -p & 0\\
0 & 0 & 0 & -1
\end{pmatrix}
$$
where $p=5$. Here we take $d = -1$ as before. The quadratic form satisfies the
conditions of \ref{prop_discr}, therefore the maximal $v$-discs are disjoint.
The result is shown in Figure \ref{fig_4}.

\newpage

\begin{figure}[h]
\centering
\includegraphics[width=10cm]{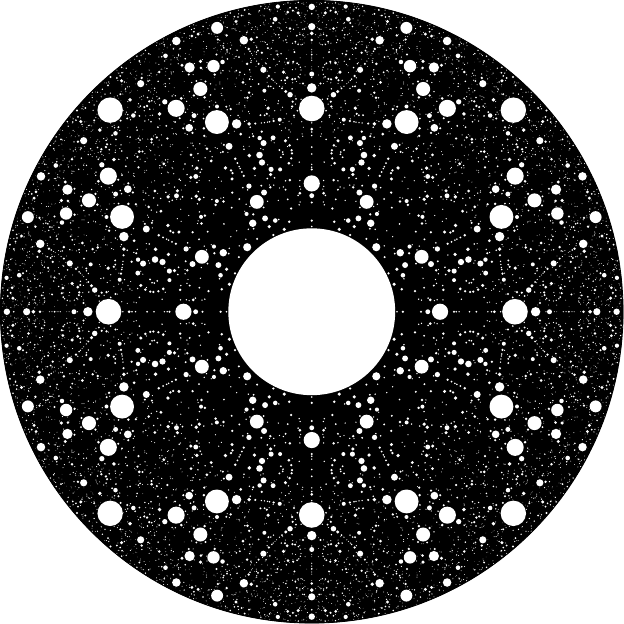}
\caption{An example where all closures of the maximal $v$-discs are disjoint (Example 4)}\label{fig_5}
\end{figure}

\subsubsection{Example 4.}

We take $d = -1$ and

$$
A = \begin{pmatrix}
3 & 0 & 0 & 0\\
0 & -1 & 0 & 0\\
0 & 0 & -63 & 0\\
0 & 0 & 0 & -63
\end{pmatrix}
$$

Note that the above quadratic form does not represent zero, as we have discussed in
\ref{example_lat}. Moreover, it satisfies the conditions of \ref{prop_discr} (taking $p=3$), and
it follows that the closures of the maximal $v$-discs are disjoint. Note that due to
the large discriminant of $q$ the vectors $v\in\Lambda$ with $q(v) = -1$ are ``rare'',
so that in the picture most of the $v$-discs have very small radii. They are, nevertheless,
everywhere dense in the plane. The big $v$-disc in the centre of Figure \ref{fig_5} is the $D_v$
with $v=(1,-2,0,0)$.

\newpage

\begin{figure}[h]
\centering
\includegraphics[width=10cm]{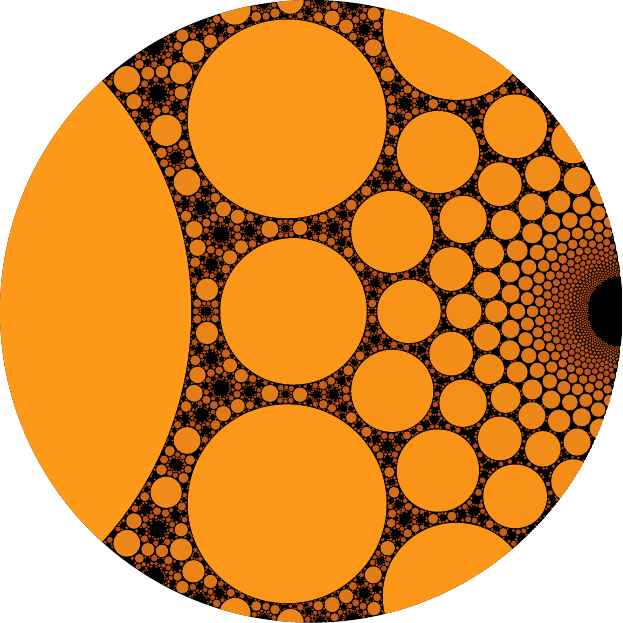}
\caption{An example where all closures of the maximal $v$-discs are disjoint (Example 5)}\label{fig_6}
\end{figure}

\subsubsection{Example 5.}

We take $d = -2$ and

$$
A = \begin{pmatrix}
-2 & 5 & 0 & 0\\
5 & 0 & 0 & 0\\
0 & 0 & -10 & 5\\
0 & 0 & 5 & -10
\end{pmatrix}
$$
The result appears in Figure \ref{fig_6}. The conditions of \ref{prop_discr} are satisfied for $p=5$.
One may guess from the picture that in this case the closures
of the $v$-discs are also disjoint, as in the previous example.

\newpage

\begin{figure}[h]
\centering
\includegraphics[width=10cm]{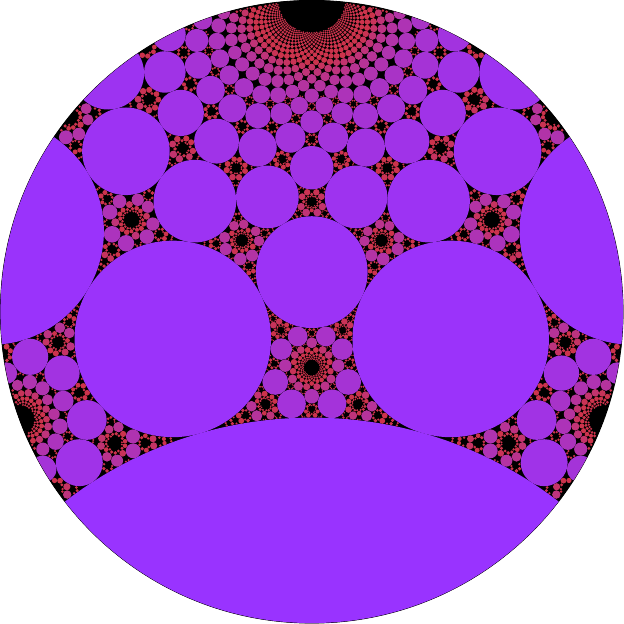}
\caption{One more example (Example 6)}\label{fig_7}
\end{figure}

\subsubsection{Example 6.}

We take $d = -1$ and

$$
A = \begin{pmatrix}
-1 & 2 & 0 & 0\\
2 & 0 & 0 & 0\\
0 & 0 & -2 & 0\\
0 & 0 & 0 & -2
\end{pmatrix}
$$
The conditions of \ref{prop_discr} are satisfied for $p=2$. The result is in Figure \ref{fig_7}.
We note that this example also appears in \cite{Su}.

\hfill

{\bf Acknowledgements:} We are grateful to Misha Belolipetsky and Misha Kapovich for useful discussions.
We are most grateful to the anonymous referees for insightful remarks and useful suggestions.

\hfill

\appendix

\section{Symmetry groups of the examples}\label{app_symmetries}

The descriptions of the symmetry groups in this appendix were kindly provided by an anonymous
referee, whom we thank for this and for many other valuable suggestions.
They describe the groups of symmetries of the configurations of
spheres in Examples 1, 2 and 5 of Subsection \ref{sec_illustr}.

All vectors below are written in the basis of $\Lambda\simeq\Z^4$ in which
$q$ is given by the matrix $A$ of the corresponding example.
For $n\in\Lambda$ with $q(n)<0$ we denote by $R_n$ the reflection in the
hyperplane $n^\perp$,
$$
R_n(x) = x - 2\frac{q(x,n)}{q(n)}\, n .
$$
For two isotropic vectors $a,b\in\Lambda$ we denote by $\mathrm{Rot}_{a,b}$
the rotation by $\pi$ about the geodesic with endpoints $[a],[b]\in\Abs$; it acts
as the identity on $\langle a,b\rangle$ and as $-\mathrm{Id}$ on $\langle a,b\rangle^\perp$.
All the elements listed below are integral isometries of $(\Lambda,q)$.

Following the referee, we say that a cusp point $p$ of a discrete group
$\Gamma\subset SO^+(1,n)$ is {\bf $k$-dimensional} if, in the upper half-space
model with $p$ at infinity, the $\Gamma_p$-orbit of an ideal point different
from $p$ generates a lattice of rank $k$ in $\R^{n-1}=\Abs\backslash\{p\}$;
here $\Gamma_p$ is the stabilizer of $p$ in $\Gamma$.
In Examples 1, 5 and 6 the ``black holes'' in the pictures correspond to
$2$-dimensional cusps, and the symmetries are easiest to see when one of
these cusps is chosen as the point at infinity.

\hfill

{\bf Example 1.} Here the group of symmetries is not generated by reflections. Let
$$
a = (1,1,0,-1),\quad b = (1,3,0,-1),\quad c = (1,1,0,1).
$$
Then the group of symmetries is
$$
\Gamma = \langle R_{(0,0,1,0)},\ R_{(0,0,0,1)},\ R_{(-1,-1,1,1)},\
\mathrm{Rot}_{a,b},\ \mathrm{Rot}_{a,c}\rangle .
$$
The point $[a]$ is a $2$-dimensional cusp.

\hfill

{\bf Example 2} (the configuration of \ref{intersect}).
Here $q=5x^2-y^2-5z^2-5w^2$, so $A$ is the diagonal matrix with entries
$5,-1,-5,-5$; the lattice of \ref{intersect} is $(\Lambda,2q)$, and rescaling
the form changes neither its isometry group nor the reflections $R_n$.
The group of symmetries is generated by the reflections $R_n$ for
$$
n\in\{(0,0,1,0),\ (0,0,0,1),\ (1,0,1,1),\ (0,0,-1,1),\ (2,5,1,0)\}
$$
and the rotation
$$
\mathrm{Rot}_{a,b}(x) = \frac{2q(a,x)b + 2q(b,x)a}{q(a, b)} - x
$$
with $a = (1,0,0,1)$, $b = (5,10,2,1)$.
The point $[a]$ is a cusp point, which is convenient to take as the point at infinity.

\hfill

{\bf Example 5.} The group of symmetries is generated by reflections, namely
by $R_n$ for
$$
n\in\{(0,0,0,1),\ (0,0,1,2),\ (0,1,1,0),\ (10,4,5,0),\ (15,6,8,4)\}.
$$
The point $(0,1,0,0)$ is a $2$-dimensional cusp: the mirrors of the first three
reflections pass through it and meet pairwise at the angles $\pi/6$, $\pi/3$ and
$\pi/2$, so its stabilizer contains a Euclidean triangle group of type $(2,3,6)$.

{

}
{\small
\noindent {\sc Ekaterina Amerik\\
{\sc Laboratory of Algebraic Geometry,\\
National Research University HSE,\\
Department of Mathematics, 6 Usacheva Str. Moscow, Russia,}\\
\tt  Ekaterina.Amerik@gmail.com}, also: \\
{\sc Universit\'e Paris-11,\\
Laboratoire de Math\'ematiques,\\
Campus d'Orsay, B\^atiment 425, 91405 Orsay, France}

\hfill

\noindent {\sc Andrey Soldatenkov\\
{\tt aosoldatenkov@gmail.com}\\
{\sc State University of Campinas \\ Rua Sérgio Buarque de Holanda, 651\\
13083-859, Campinas - SP, Brazil }}

\hfill

\noindent {\sc Misha Verbitsky\\
{\tt  verbit2000@gmail.com }}


\small
\begin{thebibliography}{MSW}

\bibitem[ALR]{ALR}
A.\ Adem, J.\ Leida, Y.\ Ruan, {\em Orbifolds and stringy topology},
Cambridge Tracts in Math., 171, Cambridge University Press, Cambridge, 2007

\bibitem[AV1]{_AV:MBM_} 
E. Amerik, M. Verbitsky, 
{\em Rational curves 
on hyperk\"ahler manifolds}, 
Int. Math. Res. Not. IMRN 2015, no. 23, 13009--13045

\bibitem[AV2]{_AV:Kaw_Mor_} 
E. Amerik, M. Verbitsky, {\em Morrison-Kawamata cone conjecture
for hyperk\"ahler manifolds,} Ann. Sci. ENS (4) 50 (2017),
no. 4, 973--993

\bibitem[AV3]{_AV:automorphisms_} 
E. Amerik, M. Verbitsky, 
{\em Construction of automorphisms of hyperk\"ahler
  manifolds}, Compos. Math. 153 (2017), no. 8, 1610--1621

\bibitem[AV4]{_AV:hyperb_}  
E. Amerik, M. Verbitsky, {\em Hyperbolic geometry of the ample
cone of a hyperk\"ahler manifold}, Res. Math. Sci. 3
(2016), Paper No. 7

\bibitem[AV5]{_AV:survey_}  
E. Amerik, M. Verbitsky, {\em Rational curves and MBM classes on hyperk\"ahler
manifolds: a survey}, in: Rationality of varieties, Progr. Math., 342, Birkh\"auser/Springer, Cham, 2021, 75--96

\bibitem[Ba1]{_Baragar1_}
A.\ Baragar, {\em The ample cone for a {$K3$} surface},
Canad. J. Math. 63 (2011) no. 3, 481--499 

\bibitem[Ba2]{_Baragar:higher_}
A.\ Baragar, {\em Higher dimensional {A}pollonian packings, revisited},
Geom. Dedicata, 195 (2018), 137--161

\bibitem[Ba3]{_Baragar:Apollonian_}
A.\ Baragar, {\em 
The Apollonian circle packing and ample cones for K3 surfaces},
arXiv:1708.06061

\bibitem[Ba4]{_Baragar:7-8_}
A.\ Baragar, {\em 
Apollonian packings in seven and eight dimensions,}
Aequationes Math. 96 (2022), no. 1, 147--165

\bibitem[Ba5]{_Baragar:Enriquez_}
A.\ Baragar, {\em 
Enriques surfaces and an Apollonian packing in eight dimensions},
Glasg. Math. J. 65 (2023), no. 1, 205--221


\bibitem[Bea]{_Beauville_} 
A. Beauville, {\em 
Vari\'et\'es K\"ahleriennes dont la premi\`ere classe de Chern est
nulle},  \textit{J. Diff. Geom.}, {\bf 18} (1983), no. 4, 755--782





\bibitem[BJ]{_Borel_Ji_}
A.\ Borel, L.\ Ji, 
{\em Compactifications of Symmetric and Locally Symmetric
  Spaces}, Mathematics: Theory \& Applications, Birkh\"auser, Boston, 2006

\bibitem[Bou]{Bou}
S.\ Boucksom, {\it Divisorial Zariski decompositions on compact complex manifolds},
Ann. Sci. \'Ecole Norm. Sup. (4) 37 (2004), no. 1, 45--76


\bibitem[CNS]{_Cano_Mavarrete_Seade_}
A.\ Cano, J.P.\ Navarrete, J.\ Seade,
{\em Complex Kleinian Groups}, Progress in Mathematics,
303,  Birkh\"auser/Springer, Basel, 2013

\bibitem[C]{_Catanese:moduli_}
F. Catanese, 
{\em A Superficial Working Guide to Deformations and Moduli},
Handbook of moduli, Vol. I, 161--215, Adv. Lect. Math. (ALM), 24,
Int. Press, Somerville, MA, 2013; arXiv:1106.1368

\bibitem[De]{_Denisi_}
F.\ Denisi, {\it Pseudo-effective classes on projective irreducible holomorphic symplectic manifolds}, Ann. Inst. Fourier (Grenoble) 2026, Online, 33 p.

\bibitem[Do]{_Dolgachev_}
I. Dolgachev, 
{\em Orbital counting of curves on algebraic surfaces and sphere packings},
 K3 surfaces and their moduli, 17--53,
Progr. Math., 315, Birkh\"auser/Springer, 2016

\bibitem[F]{_Fujiki:HK_}  
A. Fujiki, {\em On the de Rham Cohomology Group of a Compact 
K\"ahler Symplectic Manifold}, \textit{Adv. Stud.
Pure Math.}, \textbf{10} (1987), 105--165


\bibitem[Hu1]{Hu}
D.\ Huybrechts, {\em Compact hyperk\"ahler manifolds: basic results}, Invent. Math. 135
(1999), 63--113. Erratum in: Invent. Math. 152 (2003), 209--212

\bibitem[Hu2]{HuCone}
D.\ Huybrechts, {\em The K\"ahler cone of a compact hyperk\"ahler manifold},
Math. Ann. 326 (2003), no. 3, 499--513

\bibitem[Ka1]{_Kapovich:hyperbolic_}
M.\ Kapovich,
{\em Hyperbolic Manifolds and Discrete Groups,}
Modern Birkh\"auser Classics, Birkh\"auser, Boston, 2009
(reprint of the 2001 edition, Progress in Mathematics, 183)

\bibitem[Ka2]{_Kapovich:Kleinian_}
M. Kapovich, {\em Kleinian groups in higher dimensions}, in:
``Geometry and Dynamics of Groups and Spaces. In memory of
Alexander Reznikov'', M. Kapranov et al. (eds.), Birkh\"auser,
Progress in Mathematics, Vol. 265, 2007, 487--564,
available at 
{\small \url{http://www.math.ucdavis.edu/~kapovich/EPR/klein.pdf}}


\bibitem[KK]{_Kapovich_Kontorovich_}
M.\ Kapovich, A.\ Kontorovich, {\em On superintegral Kleinian sphere packings,
bugs, and arithmetic groups}, J. Reine Angew. Math. 798 (2023), 105--142

\bibitem[KN]{_Kontorovich_Nakamura_}
A.\ Kontorovich, K.\ Nakamura, {\em Geometry and arithmetic of crystallographic
sphere packings}, Proc. Natl. Acad. Sci. USA 116, no. 2 (2019), 436--441

\bibitem[Ko1]{Ko1}
S.\ Kov\'acs, {\it The cone of curves of a K3 surface}, Math. Ann. 300 (1993), no. 4, 681--691

\bibitem[Ko2]{Ko2}
S.\ Kov\'acs, {\it The cone of curves of K3 surfaces revisited}, Birational geometry, rational curves, and arithmetic, 163--169, Simons Symp., Springer, Cham, 2013


\bibitem[Ma]{_Markman:survey_}
E.\ Markman, {\em
A survey of Torelli and monodromy results for 
holomorphic-symplectic varieties}, 
Proceedings of the conference ``Complex and Differential 
Geometry'', Springer Proceedings in Mathematics, 2011, Volume 8, 257--322

\bibitem[Mat]{Mat}
D.\ Matsushita, {\it On subgroups of an automorphism group of an irreducible symplectic manifold}, arXiv:1808.10070

\bibitem[MP]{MP}
I.\ Moerdijk, D.A.\ Pronk, {\em Orbifolds, sheaves and groupoids},
K-Theory 12 (1997), no. 1, 3--21

\bibitem[Mo]{_Morris:Ratner_}
D.W.\ Morris,  
{\em Ratner's Theorems on Unipotent Flows,} 
Chicago Lectures in Mathematics, University 
of Chicago Press, 2005

\bibitem[MS]{_Mozes_Shah_} S. Mozes, N. Shah, {\em On the space of 
ergodic invariant measures of unipotent flows}, Ergodic
     Theory Dynam. Systems 15 (1995), no. 1, 149--159

\bibitem[MSW]{_Indra:MSW_} 
D.\ Mumford, C.\ Series, D.\ Wright,
{\em Indra's pearls: the vision of Felix Klein},
Cambridge University Press, 2002

\bibitem[Ni]{_Nikulin_}
V. Nikulin, {\em Integral symmetric bilinear forms and some of their applications},
Math. USSR Izv. 14 (1980), 103--167

\bibitem[Ra]{_Ratner:Raghunathan_}
M.\ Ratner, {\em Raghunathan's topological conjecture and distributions
of unipotent flows}, Duke Math. J. 63, no. 1 (1991), 235--280

\bibitem[Sa]{_Satija_}
I. Satija, {\em A tale of two fractals: 
The Hofstadter butterfly and the integral Apollonian gaskets},
 The European Physical Journal Special Topics,
 Volume 225 (2016), no. 13, pages 2533--2547


\bibitem[Sh1]{_Shah:closures_}
N.\ Shah, {\em Closures of totally geodesic immersions in manifolds
of constant negative curvature}, in: Group theory from a geometrical
viewpoint (Trieste, 1990), World Sci. Publ., River Edge, NJ, 1991, 718--732

\bibitem[Sh2]{_Shah:limiting_}
N.\ Shah,
{\em 
Limiting distributions of curves under geodesic flow on hyperbolic manifolds},
 	Duke Math. J. 148, no. 2 (2009), 251--279


\bibitem[SSV]{SSV}
N.\ Sibony, A.\ Soldatenkov, M.\ Verbitsky, {\em Rigid currents on compact hyperk\"ahler manifolds}, Ann. Sci. \'Ec. Norm. Sup\'er. (4) 58 (2025), no. 5, 1125--1154


\bibitem[So1]{So}
A.\ Soldatenkov, {\it On the Hodge structures of compact hyperk\"ahler manifolds},
Math. Res. Lett. 28 (2021), no. 2, 623--635

\bibitem[So2]{Source}
A.\ Soldatenkov, {\it Hyperbolic: a Python script that produces images
of sphere packings on the boundary of a hyperbolic space}, \url{https://github.com/aosoldatenkov/hyperbolic/}

\bibitem[St]{St}
K.\ Stange, {\it The Apollonian structure of Bianchi groups},
Trans. Amer. Math. Soc. 370 (2018), 6169--6219

\bibitem[Su]{Su}
D.\ Sullivan, {\it Entropy, Hausdorff measures old and new, and limit
sets of geometrically finite Kleinian groups},
Acta Math. 153 (1984), 259--277

\bibitem[Th]{_Thurston:thick-thin_}
W. Thurston, {\em Geometry and topology of 3-manifolds,} 1980,
Princeton lecture notes, \url{https://library.slmath.org/nonmsri/gt3m/}

\bibitem[V1]{_Verbitsky:Torelli_} 
M. Verbitsky,
{\em Mapping class group and a
  global Torelli theorem for hyperk\"ahler manifolds,} with
  an appendix by E. Markman, Duke Math. J. 162 (2013),
  no. 15, 2929--2986

\bibitem[V2]{_Verbitsky:ergodic_Erratum_} 
M. Verbitsky, {\em Ergodic complex structures on hyperk\"ahler
manifolds: an erratum}, arXiv:1708.05802

\bibitem[VW]{VW}
P.\ Vytnova, C.\ Wormell, {\it Hausdorff dimension of the Apollonian gasket},
Invent. Math. 239 (2025), no. 3, 909--946






\end{thebibliography}
\end{document}